\documentclass[12pt]{article}

\usepackage{amsfonts,epsfig,amsmath,amsthm,amssymb,graphics,verbatim,overpic,color}
\usepackage[letterpaper,margin=0.85in]{geometry}

\def\R{\mathbb{R}}
\def\C{\mathbb{C}}
\def\N{\mathbb{N}}
\def\E{\mathbb{E}}
\def\P{\mathbb{P}}
\def\I{\infty}

\newcommand{\be}{\begin{equation}}
\newcommand{\ee}{\end{equation}}
\newcommand{\bea}{\begin{eqnarray}}
\newcommand{\eea}{\end{eqnarray}}
\newcommand{\beann}{\begin{eqnarray*}}
\newcommand{\eeann}{\end{eqnarray*}}
\newcommand{\benn}{\begin{equation*}}
\newcommand{\eenn}{\end{equation*}}

\def\ra{\rightarrow}
\def\I{\infty}
\def\Id{{\textnormal{Id}}}
\def\fs{{\mathfrak{s}}}

\newcommand{\cC}{{\mathcal C}}  % calligraphic C
\newcommand{\cD}{{\mathcal D}}  % calligraphic D
\newcommand{\cF}{{\mathcal F}}  % calligraphic F
\newcommand{\cH}{{\mathcal H}}  % calligraphic H
\newcommand{\cI}{{\mathcal I}}  % calligraphic I
\newcommand{\cK}{{\mathcal K}}  % calligraphic K
\newcommand{\cL}{{\mathcal L}}  % calligraphic L
\newcommand{\cM}{{\mathcal M}}  % calligraphic M
\newcommand{\cO}{{\mathcal O}}  % calligraphic O
\newcommand{\cS}{{\mathcal S}}  % calligraphic S
\newcommand{\cU}{{\mathcal U}}  % calligraphic U
\newcommand{\cX}{{\mathcal X}}  % calligraphic X
\newcommand{\cY}{{\mathcal Y}}  % calligraphic Y
\def\txtb{{\textnormal{b}}}
\def\txtc{{\textnormal{c}}}
\def\txtd{{\textnormal{d}}}
\def\txte{{\textnormal{e}}}
\def\txti{{\textnormal{i}}}

\def\txtD{{\textnormal{D}}}

\renewenvironment{thebibliography}[1]{
  \begin{oldthebibliography}{#1}
    \setlength{\itemsep}{0.02em}
    \setlength{\parskip}{0.02em}
}
{
  \end{oldthebibliography}
}

\begin{document}

\author{Christian Kuehn\thanks{Technical University of Munich, Faculty of Mathematics,
Boltzmannstr.~3, 85748 Garching b.~M{\"u}nchen, Germany; e-mail: ckuehn@ma.tum.de}}
 
\title{Travelling Waves in Monostable and Bistable\\ Stochastic Partial Differential Equations}

\maketitle

\begin{abstract}
In this review, we provide a concise 
summary of several important mathematical results for stochastic travelling
waves generated by monostable and bistable reaction-diffusion stochastic partial 
differential equations (SPDEs). In particular, 
this survey is intended for readers new to the topic but who have some knowledge 
in any sub-field of differential equations. The aim is to bridge different backgrounds 
and to identify the most important common principles and techniques currently applied
to the analysis of stochastic travelling wave problems. Monostable and bistable 
reaction terms are found in prototypical dissipative travelling wave problems,
which have already guided the deterministic theory. Hence, we expect that these terms are 
also crucial in the stochastic setting to understand effects and to develop techniques. 
The survey also provides an outlook, suggests some open problems, and points out connections 
to results in physics as well as to other active research directions in SPDEs.
\end{abstract}


%%%%%%%%%%%%%%%%%%%%%%%%%%%%%%%%%%%%%%%%%%%%%%%%%%%%%%%%%%%%%%%%%%%%%%%%%%%%%%%%%%%%%%%%%%%%%%
\section{Introduction}
\label{sec:intro}

We consider \emph{stochastic partial differential equations (SPDEs)} of the form
\be
\label{eq:SPDEmath}
\txtd u = [\partial_{x}^2u+f(u)]~\txtd t + g(u)~\txtd W,\qquad u(0,x)=:u_0(x),
\ee
where $u=u(t,x)\in\R$, $x\in\R$, $t\in[0,\I)$ and $f$ is a given nonlinearity; 
the stochastic process $W=W(t,x)$, the map $g$, as well as the solution concept(s)
for~\eqref{eq:SPDEmath} will be specified precisely in Section~\ref{sec:SPDE}. 
The SPDE~\eqref{eq:SPDEmath} can also be written as
\be
\label{eq:SPDEphys}
\partial_t u = \partial_x^2 u + f(u) + g(u)\xi,\qquad \xi=\xi(t,x),
~\partial_t W=\xi,~u(0,x)=:u_0(x).
\ee
We focus on the classical \emph{quadratic} and \emph{cubic} nonlinearities 
given by
\benn
f(u)=f_2(u):=u(1-u)\qquad \text{or}\qquad f(u)=f_3(u):=u(1-u)(u-a),\quad a\in(0,1/2). 
\eenn
For $g(u)\equiv 0$, the equation~\eqref{eq:SPDEmath} becomes a partial 
differential equation (PDE) of \emph{reaction-diffusion type}
\be
\label{eq:PDE}
\partial_t v = \partial_x^2 v + f(v),\quad v(0,x)=v_0(x)
\ee
where we use $v=v(t,x)$ to emphasize that we work with a deterministic PDE.
Particular cases of equation~\eqref{eq:PDE} have been studied intensively for 
almost a century. Observe 
that for $f=f_2$, there are two different homogeneous steady states $v_*= 0$ and
$v_*= 1$, while in the cubic case $f=f_3$ there is the additional
steady state $v_*= a$. Looking at perturbations via $v=v_*+\varepsilon V$, 
one obtains to leading order in $\varepsilon$, the linearized system 
\be
\label{eq:PDElinearized}
\partial_t V = [\partial_x^2+\txtD_vf(v_*)]V,\qquad V=V(t,x),
~V(0,x)=V_0(x).
\ee 
Solving~\eqref{eq:PDElinearized} explicitly, e.g., via Fourier transform, one checks 
that for $f=f_2$ the state $v_*=0$ is unstable, while $v_*= 1$
is linearly stable. Hence, one refers to $f=f_2$ as the
\emph{monostable} case. For $f=f_3$, the steady states $v_*=0$ and $v_*=1$ are linearly
stable, while $v_*=a$ is unstable; hence, this case is called \emph{bistable}.
The monostable PDE~\eqref{eq:PDE} is also referred to as 
Fisher-Kolmogorov-Petrovskii-Piscounov (FKPP) equation~\cite{Fisher,
KolmogorovPetrovskiiPiscounov}. The bistable
case is called Nagumo equation in (neuro-)biology~\cite{Nagumo}, Allen-Cahn equation in
materials science~\cite{AllenCahn}, $\phi^4$-model in quantum field 
theory~\cite{Hairer1}, Schl\"ogl model in chemistry~\cite{Schloegl}, and 
real Ginzburg-Landau equation~\cite{GinzburgLandau} in the 
context of normal forms or amplitude/modulation 
equations~\cite{KuehnBook1,SchneiderUecker}.\medskip 

\begin{figure}[htbp]
	\centering
	\begin{overpic}[width=1\textwidth]{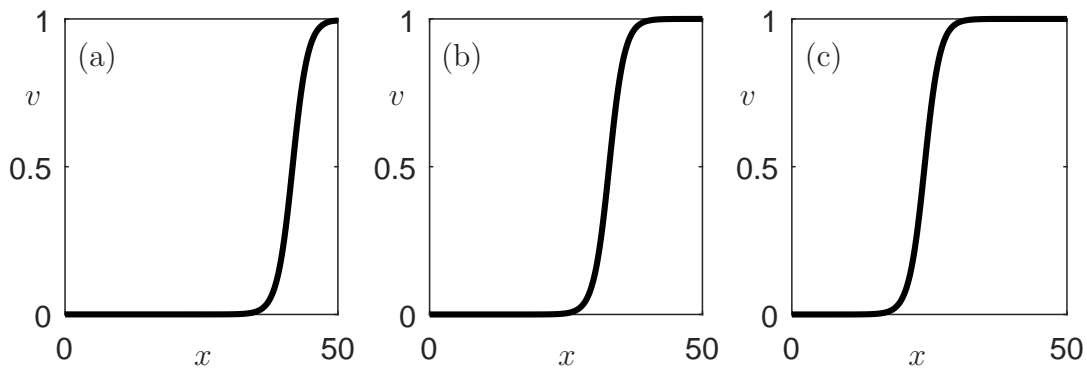}
		\put(14,23){\scalebox{1}{(a)}}
		\put(42,23){\scalebox{1}{(b)}}
	  \put(70,23){\scalebox{1}{(c)}}
		\put(23,0){\scalebox{1}{$x$}}
		\put(51,0){\scalebox{1}{$x$}}
		\put(78,0){\scalebox{1}{$x$}}
		\put(10,20){\scalebox{1}{$v$}}
		\put(38,20){\scalebox{1}{$v$}}
		\put(65,20){\scalebox{1}{$v$}}
	\end{overpic}
	\caption{\label{fig:01}Direct numerical simulation via a spatial finite-difference 
	discretization and implicit Euler time-stepping for~\eqref{eq:PDE}
	with $f=f_3$ and $a=\frac14$ on a domain $x=[0,50]$ with Neumann boundary 
	conditions. Three different times are shown: (a) $t=0$, (b) $t=25$, (c) $t=50$. We 
	clearly observe a left-travelling wave; boundary effects do not play 
	a role in practice as long as the travelling front is separated from 
	the boundary.}
\end{figure}

For the PDE~\eqref{eq:PDE} a very important class of non-steady solutions are 
\emph{travelling waves}, i.e., solutions of the form
\be
\label{eq:twave}
v(t,x)=v(x-st)=v(\eta),\qquad \eta:=x-st,
\ee
where $s\in\R$ is the \emph{wave speed}. If $s>0$ (resp.~$s<0$) then the wave
is right-moving (resp.~left-moving), while for $s=0$ we have a standing wave;
see Figure~\ref{fig:01} for an example of a travelling wave computed for the
bistable case.
 
Existence and stability of travelling
waves for the PDE~\eqref{eq:PDE} are well-studied. We recall certain parts of 
these results in Section~\ref{sec:PDE}. For the SPDE~\eqref{eq:SPDEmath} a lot 
less is known rigorously about travelling waves. However, considerable insight 
has been gained already from the perspective of physical intuition, asymptotic
approximations, direct numerical simulations, and complete proofs for certain
particular cases/aspects. The following key questions have been studied:

\begin{itemize}
 \item[(Q1)] How to define a stochastic travelling wave?
 \item[(Q2)] When is a stochastic wave 'stable'?
 \item[(Q3)] Which equation does the speed of a stochastic wave satisfy?
 \item[(Q4)] Which equation do spatial fluctuations around a stochastic wave satisfy?
 \item[(Q5)] What is the asymptotic expansion of speed and fluctuations for small noise?
 \item[(Q6)] Can we provide a bifurcation theory for stochastic waves?
\end{itemize}

As an orientation we briefly highlight already in a non-technical way, which partial
answers have been given for (Q1)-(Q6). Regarding (Q1), one can define waves via 
observables in the spatial domain, which track the wave, e.g., via point markers for 
their support or via certain norms on compact subsets of the solution. With respect to 
(Q2), one notices that the stability problem of the pattern can be posed on finite, 
as well as on infinite, time horizons. Many results show that deterministic
stability of the wave to perturbations transfers to a stable waveform also in the 
stochastically forced setting in many cases, e.g., for sufficiently small noise.
For (Q3)-(Q4), it is important to emphasize that speed and perturbations are
both random variables. The speed generally satisfies a scalar stochastic ordinary
differential equation (SODE), while the fluctuations around a reference profile
satisfy an SPDE. For (Q5), it has been proven that the asymptotics of speed 
and fluctuations differ substantially between the monostable and bistable cases. 
Large correction terms in comparison to the PDE appear for the monostable case, while
the bistable case tends to mainly produce small corrections. The question (Q6) is
least studied as phenomena such as generation of waves via bifurcations or propagation
failure of waves have only been understood in key examples. More generally, many open 
problems remain in theory of pattern formation for SPDEs. In fact, the area has been 
recognized for quite some 
time as one of the most challenging and fundamental research frontiers\footnote{To 
quote an opinion from 1991 by Glimm~\cite{Glimm}: ``The major problems for partial 
differential equations are either nonlinear or stochastic or both.'' Pattern formation
for SPDEs evidently involves the study of nonlinear stochastic PDEs.}.\medskip 

We briefly mention there are two important directions, which we do not detail
here:
\begin{itemize}
  \item[(I)] Travelling waves in ``random media'', or for random partial 
	differential equations (RPDEs), i.e., when random coefficients are introduced 
	to~\eqref{eq:PDE} as discussed e.g.~in~\cite{Xin}. 
	\item[(II)] Travelling waves for discrete/microscopic versions
of~\eqref{eq:SPDEmath}. We do not cover the RPDE case at all, yet we are going to
comment on the relation to microscopic models at relevant places. The reason is 
that the form of $g$ in~\eqref{eq:SPDEmath} is often derived from
microscopic considerations. We refer the reader for a detailed review from a physics 
perspective on discrete models to~\cite{Panja}, for micro-macro model limits 
to~\cite{Kotelenez}, and to references therein.
\end{itemize}
 
Having covered the relevant deterministic PDE background for travelling 
waves in 
Section~\ref{sec:PDE}, we introduce some SPDE basics in Section~\ref{sec:SPDE}.
Then Section~\ref{sec:waves} contains an
overview of current results on existence, speed, and stability of stochastic 
travelling waves for reaction-diffusion SPDEs with monostable ad bistable nonlinearities. 
We conclude in Section~\ref{sec:outlook} with a brief summary,
indicate connections to adjacent areas, and propose several directions for future 
work. We also alert the reader that this survey is written to communicate the
main objects, structures and effects encountered in the area. We do not attempt to
be monographic. The area is still lacking, at this point in
time, a comprehensive mathematical monograph collecting and/or developing proofs of 
all the major 'known' results. Yet, this survey aims to provide a first step towards 
this goal by organizing the very broad literature across disciplines in a more accessible 
first-reading format.

%%%%%%%%%%%%%%%%%%%%%%%%%%%%%%%%%%%%%%%%%%%%%%%%%%%%%%%%%%%%%%%%%%%%%%%%%%%%%%%%%%%%%%%%%%%%%%
\section{Deterministic Waves}
\label{sec:PDE}

In this section, we are just going to review the relevant results we need here for travelling 
waves of the deterministic PDE~\eqref{eq:PDE} with monostable and bistable 
nonlinearities. For a lot more on deterministic travelling waves, we refer the 
reader to~\cite{AronsonWeinberger,KuehnBook1,Sandstede1,VolpertVolpertVolpert} and 
references therein.\medskip  

First, let us remark that the existence and regularity theory of the 
PDE~\eqref{eq:PDE} is well-studied for a wide class of nonlinearities. 
For $f=f_3$, the highest-order polynomial term $-v^3$ provides 
dissipativity~\cite{Henry,Temam} while the linear part provides smoothing 
leading to global-in-time solutions in very regular function spaces. Due to 
smoothing~\cite[Sec.3.3,~Ex.3.7]{Henry}, one
can work with classical solutions for $t>0$ if the initial condition is taken 
sufficiently regular, which we are going to assume from now on, say taking 
continuous and bounded data $v_0\in C^0_\txtb(\R,\R)$. For $f=f_2$, an additional 
restriction of the initial condition to $v_0\geq 0$ leads to the same conclusion 
of global-in-time existence and regularity. The condition $v_0\geq 0$ is
often natural for modelling purposes of the FKPP equation and also implies
that $v(t,x)\geq 0$ for all $t\geq 0$ by using the maximum 
principle~\cite{Evans}.\medskip

There are three main classes of travelling waves we are going to consider: 
\emph{fronts}, \emph{pulses}, and \emph{wave trains}. If $v_{*,l}$ and $v_{*,r}$
are steady states for~\eqref{eq:PDE}, a travelling \emph{front} from $v_{*,l}$ to 
$v_{*,r}$ is a solution such that
\be
\label{eq:hetcond}
\lim_{\eta\ra -\I}v(\eta)=v_{*,l}\qquad \text{and}\qquad 
\lim_{\eta\ra +\I}v(\eta)=v_{*,r}.
\ee
We also refer to $v_{*,l}$ and $v_{*,r}$ as (left and right) \emph{endstates} of 
the wave. A travelling \emph{pulse} to a single steady state (or endstate) $v_*$ 
satisfies~\eqref{eq:hetcond} with $v_{l,*}=v_*=v_{r,*}$. A travelling \emph{wave 
train} is a spatially periodic pattern $v(\eta+\eta_0)=v(\eta)$ for some fixed 
$\eta_0>0$. Plugging in the travelling wave ansatz~\eqref{eq:twave} 
into the PDE~\eqref{eq:PDE} and using the chain rule yields
\be
\label{eq:secondorder}
-s\frac{\txtd v}{\txtd \eta}=\frac{\txtd^2 v}{\txtd \eta^2}+f(v),
\ee
which is a second-order ordinary differential equation (ODE). Re-writing this
ODE via $\txtd v/\txtd \eta=:w$, we get a planar first-order system
\be
\label{eq:ODE}
\begin{array}{lclcl}
\frac{\txtd v}{\txtd \eta}&=&\dot{v}&=& w,\\
\frac{\txtd w}{\txtd \eta}&=&\dot{w}&=& -sw-f(v).\\
\end{array}
\ee 
The equilibrium points of~\eqref{eq:ODE} lie on the line $\{w=0\}$ with $v_*$-values 
$v_*=0,1$ and $v_*=0,a,1$ for the quadratic and cubic nonlinearities. The 
condition~\eqref{eq:hetcond} is the defining property of a \emph{heteroclinic
orbit} in the system~\eqref{eq:ODE} from $(v_{l,*},0)$ to $(v_{r,*},0)$. Hence,
travelling fronts correspond to heteroclinic orbits, travelling pulses to 
\emph{homoclinic orbits}, and travelling wave trains to \emph{periodic orbits}
of the ODE~\eqref{eq:ODE}. 

\begin{figure}[htbp]
	\centering
	\begin{overpic}[width=1\textwidth]{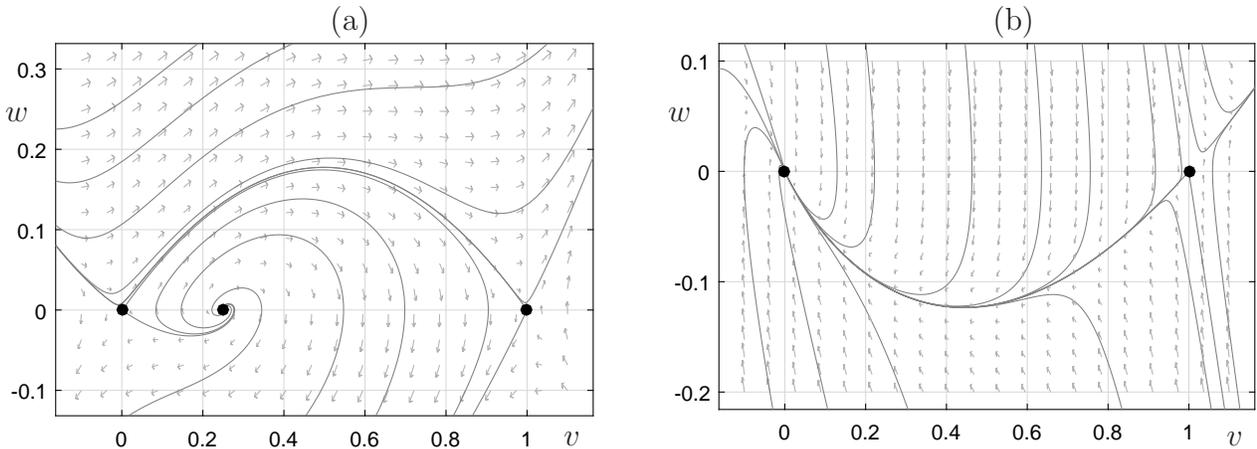}
		\put(25,32){\scalebox{1}{(a)}}
		\put(76,32){\scalebox{1}{(b)}}
		\put(43,0){\scalebox{1}{$v$}}
		\put(94,0){\scalebox{1}{$v$}}
		\put(0,25){\scalebox{1}{$w$}}
		\put(51,25){\scalebox{1}{$w$}}
	\end{overpic}
	\caption{\label{fig:02}Numerical phase portraits of the ODE~\eqref{eq:ODE}
	with steady states (black dots) and several trajectories (grey lines). (a)
	The bistable case with $f=f_3$ for $a=\frac14$ and $s=0.35$ is shown. The 
	parameter values extremely close to the existence of a heteroclinic connection 
	between the two ODE saddle equilibria at $(0,0)$ and $(1,0)$. (b) The monostable 
	case for $s=-2$, where we see a heteroclinic connection between a saddle at 
	$(1,0)$ and a node at $(0,0)$.}
\end{figure}

For $f=f_2$, one checks that $(1,0)$ is a saddle point, while $(0,0)$ is a 
stable node for $s\leq -2$, an unstable node for $s\geq 2$, a spiral sink
for $s\in(-2,0)$, a spiral source for $s\in(0,2)$, and a center for
$s=0$. Let us consider only the case $s\geq 0$ as the case $s<0$ can be dealt with
using the symmetry 
\be
\label{eq:symmetry}
(s,t,v,w)\mapsto(-s,-t,v,-w) 
\ee
of~\eqref{eq:ODE}. Although 
it is possible for $s\in[0,2)$ to construct periodic, homoclinic ($s=0$) and 
heteroclinic ($s\in(0,2)$) orbits for~\eqref{eq:ODE}, we see, due to 
the complex eigenvalues near $(0,0)$, that these orbits have negative $v$-values 
for certain $\eta$. Hence, these solutions cannot be obtained for the PDE if we 
adhere to the modelling constraint $v(0,x)=v_0(x)\geq 0$ in the monostable
case. For $s\geq 2$, one may check using~\eqref{eq:ODE} that there exists
for each fixed $s\in[2,\I)$ a unique heteroclinic orbit $\gamma_s(\eta)$ with
\benn
\lim_{\eta\ra -\I}\gamma_s(\eta)=(1,0)\qquad \text{and} \qquad 
\lim_{\eta\ra +\I}\gamma_s(\eta)=(0,0), 
\eenn 
which does satisfy the constraint $v(0,x)=v_0(x)\geq 0$; see also Figure~\ref{fig:02}(b).
This family of heteroclinic orbits $\{\gamma_s(\eta)\}_{s\geq 2}$ represents travelling 
front solutions in which 
the homogeneous state $v_*= 1$ of the original PDE \emph{invades} the homogeneous 
state $v_*= 0$. The fronts are monotone functions of $\eta$. One also refers 
to this travelling front scenario as \emph{propagation
into an unstable state} as the state $v_*= 0$ is unstable as a steady 
state of the PDE~\eqref{eq:PDE}. Of course, we may ask, which of the family of 
possible front solutions we actually observe if we consider
an initial condition for the PDE~\eqref{eq:PDE} with $f=f_2$ and $v_0\geq 0$. 
This is related to the stability question to be reviewed below.

Before discussing stability, let us also consider existence for $f=f_3$. 
For $a\in(0,1/2)$, one finds that $(v_{r,*},0)=(1,0)$ and $(v_{l,*},0)=(0,0)$ 
are saddle points, while $(v_{m,*},0)=(a,0)$ is locally stable for $s>0$, 
unstable for $s<0$, and a center for $s=0$. For $s=0$, one uses the Hamiltonian 
structure of~\eqref{eq:ODE} to show the existence of a homoclinic orbit 
representing a (standing) pulse solution. For $s\neq 0$, one can again restrict 
by the symmetry~\eqref{eq:symmetry} the range of wave speeds. Let us now take 
$s<0$. One can then prove~\cite{FifeMcLeod} that there exists for 
each $a\in(0,1/2)$ a unique $s_*=s_*(a)<0$, where~\eqref{eq:ODE} has a 
heteroclinic orbit $\gamma_a(\eta)$ such that
\be 
\lim_{\eta\ra -\I}\gamma_a(\eta)=(0,0)\qquad \text{and} \qquad 
\lim_{\eta\ra +\I}\gamma_a(\eta)=(1,0).   
\ee
The heteroclinic orbit corresponds to a front solution, which turns out to be
a monotone function of $\eta$ as a solution of the PDE. The wave speed $s_*(a)$ 
can also be expressed 
via variational principles~\cite{BenguriaDepassier,Chen1} for quite general 
bistable scalar problems. In fact, for $f=f_3$ 
one may even write explicit formulas~\cite{FifeMcLeod,Sattinger1} yielding 
\benn
\Phi(\eta)=\frac12\left[1+\textnormal{tanh}\left(\frac{\sqrt{2}}{4}\eta\right)\right],
\qquad s_*(a)=\sqrt{2}\left(a-\frac12\right),~~a\in(0,1/2).
\eenn
Yet, it is rarely a good idea if one wants to generalize arguments to rely
on these explicit formulas. Note that one can also consider the case $a\in(1/2,1)$ 
using a further 
symmetry, in which case the wave would be moving right instead of left.
Clearly, although $v_*\in \{0,1\}$ are both locally stable for the PDE~\eqref{eq:PDE} 
for $a\in(0,1/2)\cup (1/2,1)$, one stable state invades the other stable state. 
This is easy to understand using the gradient flow formulation
\benn
\partial_t v = -\nabla_{L^2(\R)}\cF(u),\qquad 
\cF(u):=\int_\R -F(u)+\frac12|\partial_x u|^2~\txtd x
\eenn  
where $F'=f$ is an anti-derivative of $f$. If $a\in(0,1/2)$, the state 
$v_*=1$ invades $v_*= 0$ as it also is the unique global 
minimum of the potential $\cF$, while the situation is reversed for
$a\in(1/2,1)$. The balanced potential case $a=1/2$ is special and leads 
to metastability if the diffusion constant in front of the Laplacian is 
small~\cite{CarrPego,KuehnBook1}. We shall not cover this metastable case 
here but see Section~\ref{sec:outlook} for further references.\medskip

The next natural question regarding the PDE~\eqref{eq:PDE} is to consider
stability of travelling waves. Obviously, non-trivial waves cannot be
globally stable in any reasonable sense since we already have at least one 
locally stable homogeneous steady state. Yet, local stability of waves, 
potentially with quite large basins of attraction, is possible. The first step 
is to consider \emph{linear stability}. Let $\Phi=\Phi(\eta)$ be a travelling 
wave, so that using the perturbation ansatz $v(t,x)=\Phi(\eta)+\varepsilon V(t,\eta)$
we get the linearized problem
\be
\label{eq:linPDEtw}
\partial_t V = \partial_\eta^2 V + s\partial_\eta V+ f'(\Phi(\eta))V=:\cL V.
\ee
We recall that the spectrum $\sigma(\cL)$ of a linear operator 
$\cL:\cX\ra \cY$, where $\cX,\cY$ are suitable Banach or Hilbert spaces, 
consists of all $\lambda\in \cC$ such that $(\lambda \Id-\cL)^{-1}$ has
no bounded inverse. We can decompose the spectrum $\sigma(\cL)=
\sigma_{\textnormal{pt}}(\cL)\cup \sigma_{\textnormal{ess}}(\cL)$, where the
\emph{essential spectrum} $\sigma_{\textnormal{ess}}(\cL)$ denotes all
$\lambda\in\C$ such that $\cL:\cX\ra \cY$ is not a Fredholm operator\footnote{There
are several slightly different definitions of the essential spectrum of an
operator, so one should check carefully, which definition each author uses.}.
The \emph{point spectrum} $\sigma_{\textnormal{pt}}(\cL)=\sigma(\cL)\setminus 
\sigma_{\textnormal{ess}}(\cL)$ consists of all eigenvalues $\lambda$ of finite 
multiplicity solving the eigenvalue problem
\be
\label{eq:Sturm}
\cL V= \lambda V,\qquad \lambda\in \C,~V\in\cX,
\ee
For our setting, one observes that $\cL$ is a special case
of a Sturm-Liouville differential operator
\be
\cL=\partial_\eta^2 + a_1(\eta) \partial_\eta +a_0(\eta)
\ee
with $a_1(\eta)= s$ and $a_0(\eta)= f'(\Phi(\eta))$, which are 
both coefficients which decay exponentially to asymptotic values 
if $\eta\ra \pm\I$. Furthermore, we
note that the travelling wave $\partial_\eta \Phi$ is always an eigenfunction with 
eigenvalue $\lambda=0$ since we may use that a travelling wave solves the PDE
and differentiate to obtain
\be
0=\partial_\eta^2 \Phi + s\partial_\eta\Phi+f(\Phi) \qquad \Rightarrow \quad 
0=\partial_\eta^2 (\partial_\eta \Phi)+s\partial_\eta(\partial_\eta \Phi)
+f'(\Phi)\partial_\eta\Phi.
\ee
The direction associated to $\partial_\eta \Phi$ corresponds to the neutral
direction induced by \emph{translation symmetry}, i.e., if $\Phi(x-st)$ is a 
travelling wave then so is $\Phi(x-st+\eta_0)$ for $\eta_0\in\R$ fixed. The
neutral mode does not contribute to the linear stability analysis for the
PDE~\eqref{eq:PDE}; it is also referred to as the \emph{Goldstone mode} in
the physics literature. If we remove the eigenvalue associated to the Goldstone 
mode from the spectrum $\sigma(\cL)$ and if we can prove the remaining part of the 
spectrum is contained properly in the left half of the complex plane, then the 
wave is called \emph{linearly stable}. Linear stability is a necessary condition 
to obtain nonlinear stability, i.e., that that there exists a constant $\eta_0>0$ 
such that
\be
\label{eq:typstabres}
\lim_{t\ra +\I}\| v(t,\cdot)-\Phi(\cdot-st+\eta_0)\|=0
\ee
where  $\|\cdot\|$ is a norm on the (spatial) 
function space, and the initial condition $v_0$ of the solution $v$ is taken 
from a certain class of data within a basin of attraction of the wave. 
So stability means that initial conditions in some set around the wave converge 
to \emph{some translate} of the wave. Returning to linear stability,
the first standard setting is to consider $\cX=H^2(\R)$ and $\cY=L^2(\R)$. 
In fact, the \emph{essential spectrum} for our case can be inferred 
from the \emph{asymptotic linearized operators}
\be
\label{eq:twasymp}
\cL_\pm:=\partial_\eta^2 + s \partial_\eta +f'(\Phi(\pm\I)).
\ee
To see this in a bit more detail, consider the linear problems
\be
\label{eq:linPDE1}
\partial_t V= \cL_\pm V 
\ee
We use the ansatz $V(t,\eta)=\txte^{\txti k \eta - \lambda t}$ (where 
$\txti:=\sqrt{-1}$) in the linear PDE~\eqref{eq:linPDE1} and use the chain rule. 
Upon dividing each side by $\txte^{\txti k \eta - \lambda t}$ one obtains the 
dispersion relations
\be
d_\pm(\txti k,\lambda):= \lambda + (\txti k)^2 + c \txti k + f'(\Phi(\pm\I)) =0.
\ee
A \emph{dispersion relation} connects the temporal decay $\lambda$ to the spatial
\emph{wave number} $k$~\cite{KuehnBook1,SchneiderUecker}. The key objects enclosing 
the essential spectrum are the parabolic curves
\be
\label{eq:disprel1}
\sigma_\pm:=\{\lambda\in \C|d_\pm(\txti k,\lambda)=0,k\in\R\}=
\{\lambda\in \C:\textnormal{Re}(\lambda) = f'(\Phi(\pm\I))-
(\textnormal{Im}\lambda/s)^2 \}.
\ee
but this is non-trivial to prove~\cite{KapitulaPromislow,Sandstede1}.
Although up to this point, the (linear) stability 
problem~\cite{Sandstede1} can be set up in the same way for $f=f_2$ and 
$f=f_3$, the actual spectra turn out to differ substantially. 

For $f=f_3$, one can prove that $\sigma_{\textnormal{ess}}(\cL)$ is properly
contained in the left-half of the complex plane for the pulse and front 
solutions if $s\leq 0$. Indeed, recall that $f=f_3$, we have $f'(\Phi(\pm\I))<0$ as 
both endstates are stable. Therefore, the curves~\eqref{eq:disprel1} both
lie in the left half of the complex plane. So spectral stability is completely 
determined by the eigenvalue problem~\eqref{eq:Sturm}. For the travelling front 
and the travelling pulse, it is possible to check that the point spectrum consists 
of finitely many eigenvalues using Sturm-Liouville 
theory~\cite{KapitulaPromislow}. The standing pulse for $s=0$ is unstable, 
while the travelling front for $s_*=s_*(a)<0$ is linearly/spectrally stable 
see e.g.~\cite[Sec.9]{KuehnBook1}. In fact, there is a \emph{spectral gap} 
for the pulse in the bistable case, i.e., there exists 
a fixed constant $C_*=C_*(a)>0$ such that if $\lambda\in\sigma(\cL)\setminus\{0\}$
then $\textnormal{Re}(\lambda)<-C_*$. From this spectral gap, and upon 
factoring out the translation symmetry direction, one can even establish \emph{nonlinear
orbital asymptotic stability} of the travelling pulse~\cite{KapitulaPromislow,Henry}, 
i.e., for the full nonlinear PDE we have that a small perturbation to a travelling 
wave converges as $t\ra +\I$ to a translate of the travelling wave; see 
also~\cite[Def.~4.3.4]{KapitulaPromislow} and equation~\eqref{eq:typstabres}.

For $f=f_2$, the stability question is very different. Clearly, the essential
spectrum $\sigma_{\textnormal{ess}}(\cL)$ now also contains parts in the right
half of the complex plane as for the unstable endstate we have 
\be
\sigma_-=
\{\lambda\in \C:\textnormal{Re}(\lambda) = 1-(\textnormal{Im}\lambda/s)^2 \}.
\ee
However, one may use suitable exponentially weighted spaces $\cX_w$ and $\cY_w$ 
to shift $\sigma_{\textnormal{ess}}(\cL)$ for $\cL:\cX_w\ra \cY_w$ to a half-plane 
$\{\lambda\in\C:\textnormal{Re}(\lambda)<-C\}$ for some $C>0$. As an example, we 
can take $\cX_w=H^2(\R;w)$ with norm
\benn
 \|V\|_{H^2(\R;w)}=\left\|V(\cdot)(1+\txte^{-c\cdot})\right\|_{H^2(\R)},\qquad 
V=V(x),
\eenn
and we take $w(x)=\txte^{-cx}$ for some suitable constant $c>0$. Similarly, one
can define $\cY_w=L^2(\R;w)$. Hence, stability depends again on the eigenvalue 
problem $\cL V=\lambda V$. One may prove that if the wave front $\Phi$ decays faster than
$\txte^{-cx}$ as $x\ra \I$, then it is asymptotically stable (up to shifts) in the
norm with weight $w(x)=\txte^{-cx}$. Furthermore, for $s>0$ and $v_0\geq 0$, 
only the travelling front with \emph{minimal speed} $s=2$ turns out to be locally 
stable for all initial data, which have at least exponential tails 
\be
|v(0,x)-\Phi(x+\eta_0)|=\cO(\txte^{-x})\quad \text{as $x\ra \pm\I$}
\ee
However, one may select other initial data to achieve different/faster 
speeds. For certain compactly supported initial conditions sufficiently close 
to front-like profiles at the boundary of their support, one can prove that 
there is a left-moving and a right-moving front, 
each with minimal wave speed $|s|=2$. In particular, front solutions with minimal 
wave speed are selected dynamically by the PDE evolution from the infinite family 
of possible FKPP front profiles.\medskip 

The front solution for the monostable case is also called a \emph{pulled
front} as its wave speed corresponds precisely to the linear spreading 
speed of small perturbations near the unstable steady state. More precisely, 
the speed of the front can be calculated from the linearization
around the unstable state as follows: consider the general linearization 
around a steady state $v_*$ is given by~\eqref{eq:PDElinearized}, which can 
be written as
\be
\label{eq:linFKPP}
\partial_t V = \partial_x^2 V + \txtD f(v_*)V =: \partial_x^2 V + a_* V.
\ee
Let $\{(t,x_\kappa(t)):t\in[0,T)\}$ be a level curve such that $V(t,x_\kappa(t))=\kappa
\in[0,1]$ and define a \emph{linear spreading speed} by
\benn
s_*:=\lim_{t\ra +\I} \frac{\txtd x_\kappa}{\txtd t}.
\eenn 
Furthermore, applying the \emph{Fourier transform} 
\benn
\hat{V}(t,k):=\int_\R \txte^{-\txti kx} V(t,x)~\txtd x
\eenn
to~\eqref{eq:linFKPP} and substituting the ansatz 
$\hat{V}(t,k)=\hat{V}(0,k)\exp(\txti \omega_*(k)t)$ into the resulting equation, yields
the \emph{dispersion relation}
\be
\label{eq:dispersionrelation}
\omega_*(k)=-\txti (k^2-1)
\ee
relating the (complex) \emph{frequency} $\omega_*$ to the \emph{wave number} $k$.
A front is pulled~\cite{vanSaarloos} if its speed coincides with a linear spreading 
speed, which can just be calculated from the dispersion relation by the relations
\be
\label{eq:Saarloos}
s_*=\frac{\txtd \omega_*}{\txtd k}(k_*)\qquad \text{and}\qquad 
s_*=\frac{\textnormal{Im}[\omega_*(k_*)]}{\textnormal{Im}[k_*]}, 
\ee
where $k_*\in\C$ is a constant also called the \emph{linear spreading point}. 
Intuitively, a linear spreading point describes the dominant wave number of the 
linearized dynamics near the unstable state. One
easily checks that for the monostable case the linearization at the unstable state 
yields $a_*=1$ and $\omega_*(k)=\txti(1-k^2)$. Therefore, the 
conditions~\eqref{eq:Saarloos} give $k_*=\pm \txti$ and $s_*=\pm 2$. Hence, we
precisely recover the minimal wave speed for the monostable case, so the front
is \emph{pulled}. For the bistable case $f=f_3$, one can carry out the same 
linear calculation but finds that the linear spreading speed does not coincide with
the true wave speed. In this case, the front is called \emph{pushed}.

%%%%%%%%%%%%%%%%%%%%%%%%%%%%%%%%%%%%%%%%%%%%%%%%%%%%%%%%%%%%%%%%%%%%%%%%%%%%%%%%%%%%%%%%%%%%%%
\section{SPDE Background}
\label{sec:SPDE}
 
To define the noise, let $\cH$ be a Hilbert space with inner product $\langle\cdot,
\cdot\rangle_\cH$ and suppose $\cH$ is also a suitable function spaces over 
$\R$, e.g., we may simply think of $\cH=L^2(\R)$ as one key example but also suitably
weighted Lebesque or Sobolev spaces are frequently used to allow for non-vanishing functions 
at $x=\pm \I$ or to lift technical results from bounded domains to unbounded ones. 
Next, we consider an operator 
$Q:\cH\ra \cH$ and assume it has eigenfunctions $\{e_k\}_{k=1}^\I$ and associated non-negative 
eigenvalues $\{\lambda_k\}_{k=1}^\I$ such that $Qe_k=\lambda_k e_k$ for each $k\in \N$. We
define~\cite{DaPratoZabczyk} an $\R$-valued \emph{$Q$-Wiener process} as
\be
\label{eq:sumnoise}
W(t,x):=\sum_{k=1}^\I \sqrt{\lambda_k}e_k(x)B_k(t),\qquad W(t):=W(t,\cdot),
\ee 
where $\{B_k\}_{k=1}^\I$ are independent identically distributed (iid) Brownian motions over
a probability space $(\Omega,\cF,\P)$. Using standard properties of Brownian motion, it is
easy to see that we have zero mean $\E[W(t)]=0$ and the \emph{correlation function} $\E[W(t,x)W(s,y)]
=\min(t,s)q(x,y)$, where the spatial correlation function $q$ is determined via $Q$ by
\benn
Q\zeta(x)=\int_\R \zeta(y)q(x,y)~\txtd y.
\eenn
Hence, $Q$ can be viewed as a \emph{covariance operator} since
\benn
\E[\langle W(t),\zeta_1\rangle_\cH \langle W(s),\zeta_2\rangle_\cH]=
\min(t,s)~ \langle \zeta_1,Q\zeta_2\rangle_\cH
\eenn
for all $\zeta_{1,2}\in\cH$. Let $g:\R\ra \R$ be a given
function and now we may consider the SPDE~\eqref{eq:SPDEmath} as 
an \emph{evolution equation} on $\cH$ for $u=u(t)=u(t,x)$ as follows
\be
\label{eq:SPDEmath1}
\txtd u = [\partial_{x}^2u+f(u)]~\txtd t + g(u)~\txtd W,\qquad u(0,x)=:u_0(x).
\ee 
In more generality, one can also use maps taking values in linear operators instead of $g$, 
i.e., one may take in the noise term $G:\cH\ra \cL(\cU,\cH)$ for some Hilbert space $\cU$ so that 
$(G(u)\zeta)(x)=g(u(x))\zeta(x)$ with $\zeta\in \cU$. However, we shall not need this more 
general viewpoint here. The SPDE~\eqref{eq:SPDEmath1} can be written formally as
\be
\label{eq:SPDEphys1}
\partial_t u = \partial_x^2 u + f(u) + g(u)\xi,\qquad \xi=\xi(t,x),
~\partial_t W=\xi,~u(0,x)=:u_0(x),
\ee
where $\xi$ is just a \emph{generalized stochastic process}~\cite{ArnoldSDEold}. Indeed, 
Brownian motion has only H\"older-regularity in time of $(\frac12-\rho)$ for any $\rho>0$, 
so its time 
derivative is a generalized function/process. There are two common assumptions on $Q$. If $Q=\Id$, 
then we have \emph{space-time white noise} with $\E[\xi(t,x)\xi(s,y)]=\delta(t-s)\delta(x-y)$, where 
$\delta$ is the usual Dirac-delta generalized function, and we have derived the correlation function 
of $\xi$ from the one for $W$ by taking time derivatives. If $Q$ is a trace-class operator, i.e., 
$\sum_{k=1}^\I\lambda_k<+\I$, then we have a \emph{spatially-correlated noise} with 
$\E[\xi(t,x)\xi(s,y)]=\delta(t-s)q(x,y)$ for some spatial correlation function, which we
are just going to assume to
depend just on the difference between spatial locations 
\benn
q(x,y)=q(x-y) 
\eenn
from now on. Spatially correlated
noise has higher regularity than space-time white noise; see also~\cite{DaPratoZabczyk,Hairer1} 
for more details. One possible solution concept to the SPDE~\eqref{eq:SPDEmath1} is to 
consider a \emph{mild solution}
\be
\label{eq:solmild}
u(t)=S(t)u_0+\int_0^t S(t-s)f(u(s))~\txtd s + \int_0^t S(t-s)g(u(s))~\txtd W(s),
\ee 
where $S(t)=\txte^{t\partial_x^2}$ is the usual semigroup generated by the Laplacian on 
$\cD(\partial_x^2)$, where
a special choice of (weighted) space is often necessary for the semigroup domain since travelling 
fronts with non-zero endstates do not even lie in standard Lebesgue spaces such as $L^2(\R)$. One 
interpretation of the last 
integral in~\eqref{eq:solmild} is as an \emph{It\^o integral}~\cite{DaPratoZabczyk}. As usual 
in the theory of stochastic differential equations~\cite{Oksendal,Protter}, we could also consider the Stratonovich 
form of~\eqref{eq:SPDEmath1}
\be
\label{eq:SPDEmath2}
\txtd u = [\partial_{x}^2u+\tilde{f}(u)]~\txtd t + g(u)\circ\txtd W,\qquad u(0,x)=:u_0(x),
\ee 
where the integral in an analogous mild solution formula~\eqref{eq:solmild} has to be interpreted
as a \emph{Stratonovich integral} if~\eqref{eq:SPDEmath2} is used and $g(u)$ can be a linear operator. 
For trace-class noise one has the formal relation~\cite{DuanWang}
\be
\label{eq:Itoformula}
\tilde{f}(u)=f(u)+\frac12\sum_{k=1}^\I\sqrt{\lambda_k}((\txtD_u g)(u)(g(u)e_k))e_k,
\ee
where $g(u)e_k\in\cH$, $\txtD_u g$ denotes the Fr\'echet derivative, and $(\txtD_u g)(u)(\cdot)$ is
again a linear operator. In summary, an \emph{It\^o-Stratonovich correction/conversion term} of the
usual form anticipated from SODEs appears. 
The rigorous derivation of results such as~\eqref{eq:Itoformula} requires an infinite-dimensional 
\emph{It\^o formula}, i.e., an infinite-dimensional stochastic version of the chain rule. In general, 
it is not easy to prove rigorous It\^o-type formulas for solutions of SPDEs; see the review in 
the introduction of~\cite{DaPratoJentzenRoeckner}. Furthermore, correction terms obtained in 
(numerical) approximations of SPDEs pose similar technical issues~\cite{HairerMaas}. It is important 
to point out that any It\^o formula, as well as any It\^o-Stratonovich correction such 
as~\eqref{eq:Itoformula}, is expected to be ``quadratic'' in the noise term $g$. Therefore, if
$g$ is multiplied by a small scalar $\sigma>0$, then we obtain terms of order $\sigma^2$. It is
hence natural to conjecture that in the small noise regime $\cO(\sigma^2)$-terms only have a 
major impact on the dynamics if there is some form of instability in the problem; see 
also~\cite{KuehnCurse} and the last part of Section~\ref{sec:outlook}.

Another frequently used solution concept for \eqref{eq:SPDEmath1} are \emph{weak} or
\emph{variational solutions}
\be
u(t)= \langle \zeta ,u_0\rangle_\cH + \int_0^t \langle \partial_x^2 \zeta,u(s)\rangle_\cH +
\langle \zeta,f(u(s))\rangle_\cH~\txtd s 
+ \int_0^t \langle\zeta, g(u(s))~\txtd W(s)\rangle_\cH,\quad
\forall \zeta\in\cD(\partial_x^2),  
\ee
where $\langle\cdot,\cdot\rangle_\cH$ again denotes the inner product in $\cH$. Under 
reasonable conditions 
on $f$, $g$, and $Q$, one can show that local-in-time mild solutions~\cite{DaPratoZabczyk} 
and weak solutions exist; under reasonable technical conditions these solutions even 
coincide in many cases~\cite{PrevotRoeckner}. 

On a technical level, the \emph{uniqueness} or the \emph{weak uniqueness} of solutions 
is often needed to ensure the \emph{strong Markov property}, which - roughly speaking - 
means that the solution $u(t+\tau)$ for $t>0$ and a stopping time $\tau$ is independent 
of $u(\tau)$ conditional upon $\tau<+\I$. For general nonlinear non-Lipschitz
reaction terms, such as $f=f_2$ and $f=f_3$, it is usually hard to prove uniqueness of solutions 
directly. Yet, first proving global-in-time existence~\cite{CrauelFlandoli} and 
then cutting off the nonlinearity if $|u|$ is large to obtain a globally Lipschitz 
problem~\cite{Chow}, provides the same work-around well-known for deterministic dissipative 
reaction-diffusion 
PDEs~\cite{KuehnBook1,Robinson1,Temam}. Of course, one should be aware that for large noise, 
rare event fluctuations, and/or finite-time blow-up scenarios, the approach of cutting off the 
nonlinearity can influence the dynamics significantly. 

There are other solution concepts such
as strong, kinetic~\cite{DebusscheHofmanovaVovelle}, martingale~\cite{BrzezniakGatarek}, 
pathwise mild~\cite{KuehnNeamtu}, or renormalized~\cite{Hairer1} solutions. Strong solutions 
rarely exist~\cite{DaPratoZabczyk} due to the roughness of the noise. The other classes 
of solutions for SPDEs are more complicated to define/construct. Furthermore, they are not 
immediately needed here to make sense of uniformly-parabolic scalar reaction-diffusion SPDEs 
on $\R$ and their travelling wave-type solutions; see Section~\ref{sec:outlook} for cases,
where other solution concepts enter the picture. 
Yet, what is often needed are \emph{comparison principles}. For 
example, if we assume that $f=f_2$ or $f=f_3$, and
\be
\label{eq:baseinv}
g(0)= 0,\qquad g(1)= 0,\qquad u_0(x)\in[0,1],
\ee
then we intuitively expect that $u(t,x)\in [0,1]$ for all $t>0$ almost surely. An even 
stronger result would be a comparison theorem, i.e., one assumes that
\benn
0\leq u_1(0,x)\leq u_2(0,x)\leq 1	
\eenn
is true and then concludes that 
\benn
0\leq u_1(t,x)\leq u_2(t,x)\leq 1,\qquad t>0,	
\eenn
also holds almost surely. These types of invariant-region and comparison results indeed hold under 
various technical conditions; see~\cite{Assing,DonatiMartinPardoux,Gyongy1,MantheyZausinger,Mueller2,
MuellerPerkins,Kotelenez1,Shiga,TessitoreZabczyk} and references therein. If $u(t,x)\in[0,1]$ for 
$t>0$ then it is easier to study the long-time asymptotic behaviour of travelling wave-type 
solutions as one already has a-priori boundedness and one can simply cut off the nonlinearity 
outside of $u\in[0,1]$ as necessary. 
However, for certain types of noises, such as additive noise given by $g\equiv \text{const.}$, 
simple comparison principles are rarely helpful as large deviations of the solution from 
the region $u\in[0,1]$ are 
are going to occur with positive probability. In this situation, it often makes more sense 
to focus on the behaviour of travelling wave-type solutions for times $t\in[0,T]$ for some 
fixed \emph{finite time scale} $T>0$.

%%%%%%%%%%%%%%%%%%%%%%%%%%%%%%%%%%%%%%%%%%%%%%%%%%%%%%%%%%%%%%%%%%%%%%%%%%%%%%%%%%%%%%%%%%%%%%
\section{Stochastic Waves}
\label{sec:waves}

\begin{figure}[htbp]
	\centering
	\begin{overpic}[width=1\textwidth]{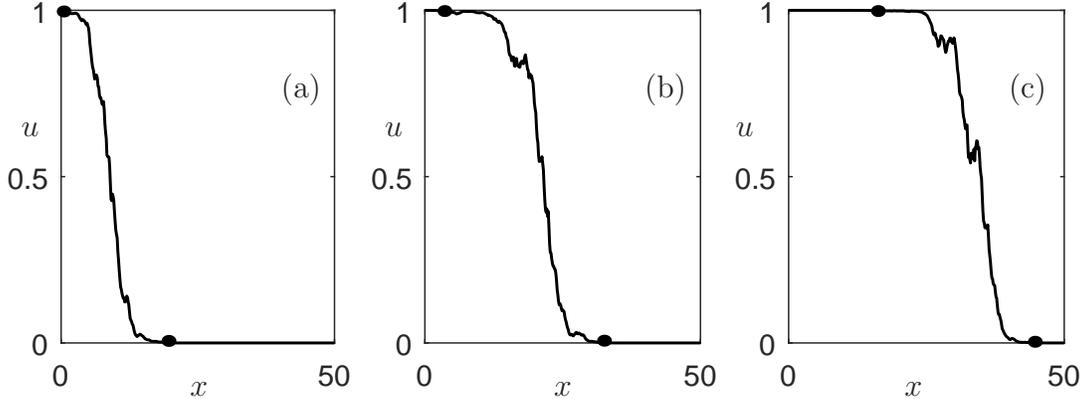}
		\put(30,23){\scalebox{1}{(a)}}
		\put(58,23){\scalebox{1}{(b)}}
	  \put(86,23){\scalebox{1}{(c)}}
		\put(23,0){\scalebox{1}{$x$}}
		\put(51,0){\scalebox{1}{$x$}}
		\put(78,0){\scalebox{1}{$x$}}
		\put(10,20){\scalebox{1}{$u$}}
		\put(38,20){\scalebox{1}{$u$}}
		\put(65,20){\scalebox{1}{$u$}}
	\end{overpic}
	\caption{\label{fig:03}Direct numerical simulation via a spatial finite-difference 
	discretization, implicit Euler time-stepping (``Milstein''), and noise for the FKPP 
	equation~\eqref{eq:SPDEFKPP1} with $g(u)=u(1-u)$ on a domain $x=[0,50]$ with 
	Neumann boundary conditions. The noise is chosen via an operator $Q$, which is 
	simultaneously diagonalizable with the Laplacian with eigenfunctions $e_k$. The noise is 
	white but truncated after $25$ eigenmodes, i.e., $\lambda_k=1$ for $k\leq 25$ and 
	$\lambda_k=0$ otherwise. The initial condition $u_0$ is chosen so that $u_0(x)=1$ for 
	$x\in [0,2.5]$ and zero otherwise. The equation satisfies an invariant region principle 
	with solutions $u\in[0,1]$. Three different times are shown: (a) $t=6.25$, (b) $t=12.5$, 
	(c) $t=18.75$. We clearly observe a noisy right-travelling stochastic wave, which actually 
	turns out to have well-defined wave speed limits~\eqref{eq:timelimits}; the markers $a(t)$
	and $b(t)$ defined in \eqref{eq:abdef} are shown as black dots. For more on
	SPDE numerics we refer to~\cite{LordPowellShardlow}.}
\end{figure}

It is already non-trivial, how to define a \emph{stochastic travelling wave} in the 
context of time-dependent stochastic forcing and roughness of the solution. Several related
approaches exist, which we briefly review. Suppose we use a deterministic continuous initial 
condition $u_0$ resembling a front such that for some $K_0>0$ we have
\be
\label{eq:intlocsupp}
u_0(x)=\left\{\begin{array}{ll} 1~ & \text{for $x<-K_0$,}\\
0~ & \text{for $x>K_0$,}\\
\end{array}\right.
\ee
and $u_0(x)\in[0,1]$ for all $x\in\R$. Suppose we select the noise term $g$ so that an 
invariant-region/comparison principle holds so $u(t,x)\in[0,1]$ for all $t\geq 0$. Furthermore,
suppose we can show that $u(-\I,t)=1$ and $u(+\I,t)=0$ for all $t\geq 0$, then the choice of 
the initial condition entails that the solution can only be different from the values of the
endstates over a bounded set for $t>0$. In
this case, it is natural to consider \emph{level sets} and define
\be
\label{eq:abdef}
a(t):=\sup \{z\in\R:u(t,x)=1,~\forall x\leq z\},\qquad b(t):=\sup \{z\in\R:u(t,x)=0,~\forall x\geq z\};
\ee
see also Figure~\ref{fig:03}.
The random variables $a(t)$ and $b(t)$ measure the spread of left and right edges of a 
travelling wave front. Of course, the convention to look at endpoints of the support is 
somewhat arbitrary and we could also consider other level sets
\benn
c_\alpha(t):= \sup \{z\in\R:u(t,x)=\alpha,~x\leq z\},\qquad \alpha\in(0,1).
\eenn
To obtain a well-defined \emph{stochastic wave speed}, one often aims to show that at least 
one of following limits exists
\be
\label{eq:timelimits}
\lim_{t\ra +\I} \frac{a(t)}{t},\qquad \lim_{t\ra +\I} \frac{b(t)}{t},\qquad
\lim_{t\ra +\I} \frac{c_\alpha(t)}{t}.
\ee
In this case, the stochastic processes $a(t)$, $b(t)$, $c_\alpha(t)$ describe the positions
so their time derivatives are stochastic processes processes, which can be taken as a
definition of the speed. We shall refer to concrete examples, how to use level sets and 
comparison principles in Section~\ref{ssec:mostoch}; see also Figure~\ref{fig:03}. 

A related, somewhat weaker notion, is to define stochastic travelling waves via 
\emph{stationary laws}~\cite{TribeWoodward}. Suppose we are interested in
the case 
\be
\label{eq:assumeTW}
W(t,x)=B_1(t),~f=f_2~\text{or}~f=f_3 \qquad \text{and}\qquad g(0)=0,~g(1)=0,\qquad u(t,x)\in[0,1], 
\ee
where the relevant solutions are fronts connecting $u=0$ to $u=1$.
Define the space
\benn
\cS:=\{\phi:\R\ra [0,1]:~\phi(-\I)=1,~\phi(+\I)=0,~\text{$\phi$ decreasing and right-continuous}\}
\eenn
with the $L^1_{\textnormal{loc}}(\R)$ topology.
Then we can define a \emph{wave marker} $C_\alpha$, similar to $c_\alpha$, for each $\phi\in\cS$ 
and center the wave accordingly
\benn
C_\alpha(\phi):= \inf\{x:\phi(x)<\alpha\},\qquad \phi_\alpha(x):= \phi(C_\alpha(\phi)+x),
\eenn
so that $\phi_\alpha$ is just $\phi$ re-centered at height $\alpha$. A stochastic travelling 
wave is a solution $u(t,x)$ for which the re-centered process
\be
\label{eq:statreq} 
u_\alpha(t,x)=u(t,x+C_\alpha(u(t))) 
\ee
is a stationary process in time. The law of the stochastic
wave is then given by the law of $u_\alpha(0,x)$. On can actually prove under the assumptions 
\eqref{eq:assumeTW} that if $u_0\in\cS$ almost surely, then $u(t)\in \cS$ almost surely.
Starting from the Heaviside function $u_0(x)=\textbf{1}_{\{x<0\}}$, assuming in addition for
$f=f_3$ that $g(a)\neq 0$, one may use a stochastic ordering technique~\cite{TribeWoodward}
to obtain stochastic travelling waves in the sense of the last definition for monostable and 
bistable cases; see also~\cite{HuangLiuWang} for an extension to $g(1)\neq 0$.\medskip

\begin{figure}[htbp]
	\centering
	\begin{overpic}[width=1\textwidth]{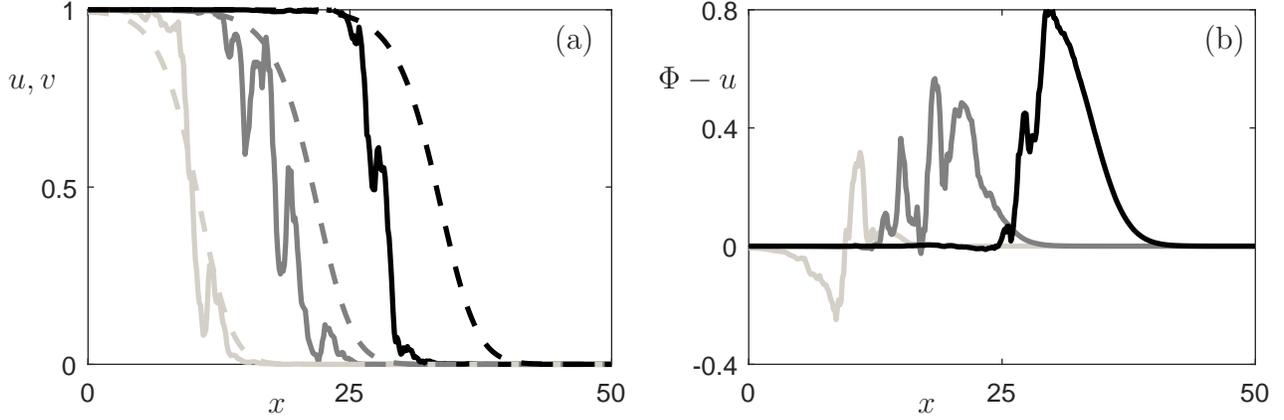}
		\put(42,28){\scalebox{1}{(a)}}
		\put(92,28){\scalebox{1}{(b)}}
		\put(20,0){\scalebox{1}{$x$}}
		\put(70,0){\scalebox{1}{$x$}}
		\put(0,25){\scalebox{1}{$u,v$}}
		\put(50,25){\scalebox{1}{$\Phi-u$}}
	\end{overpic}
	\caption{\label{fig:04}(a) Same numerical simulation as in Figure~\ref{fig:03} except that
	we have selected $\sigma=5$ for the noisy case. The three time snapshots are shown in light 
	grey, dark grey and black. The dashed curves are for $\sigma=0$ providing a deterministic
	wave $\Phi$. (b) Difference between a fixed deterministic reference wave starting from the 
	same initial data and the stochastic wave, i.e., $\Phi-u$; we
	observe that the difference can indeed get quite large in time so that an adaptive 
	minimizer as in~\eqref{eq:minLT} is useful.}
\end{figure}

Another more general, and from a deterministic viewpoint potentially more natural, approach is 
to try to define the wave and its speed via a similar strategy as for PDEs, i.e., to focus entirely 
on the moving frame dynamics. Next, we observe that stochastic and deterministic waves may 
deviate substantially if we just compare them for the same initial condition as shown in 
Figure~\ref{fig:04}(b) without re-centering. So if we want to make reasonable comparison, the 
re-centering strategy is imperative. This viewpoint requires us to analyze \emph{re-centering processes} 
in their own right without any a-priori requirements on stationarity as above for~\eqref{eq:statreq}. 
As above, since the stochastic forcing is time-dependent, we now consider a \emph{dynamic
re-centering} of the form 
\be
\label{eq:SPDEansatz}
u_p(t,x)=u(t,x-p(t)),
\ee
see also Figure~\ref{fig:04}. For the classical PDE case, we just have $p(t)=st$, where $s$ is 
a fixed \emph{wave speed} and $p(t)$ is the \emph{position of the wave}. In the 
SPDE case, $p(t)$ is generically a stochastic process. It can often be shown that
there exists an integral representation of the form
\be
\label{eq:integrated}
p(t)=\int_0^t \fs(r)~\txtd r,
\ee 
where $\fs=\fs(t)$ is another stochastic process. There is the obvious physical 
relation between the speed $\fs$ and the position $p$ in~\eqref{eq:integrated}. 
To see this in more detail, consider for simplicity the Stratonovich 
SPDE~\eqref{eq:SPDEmath2} so that we can use the standard chain rule applied 
to~\eqref{eq:SPDEansatz}. Hence, we obtain
\benn
\txtd u_p(t,x) = (\partial_x u_p)(t,x)\circ \txtd p + \txtd u(t,x-p(t)).
\eenn
A time-dependent shift does not change the covariance of $W$. Therefore, we get
in law that
\bea
\txtd u_p &=& [\partial_{x}^2u_p+\frac{\txtd p}{\txtd t}\partial_xu_p
+\tilde{f}(u_p)]~\txtd t + g(u_p)\circ\txtd W,\nonumber\\
&=&[\partial_{x}^2u_p+\fs(t)\partial_xu_p
+\tilde{f}(u_p)]~\txtd t + g(u_p)\circ\txtd W,\label{eq:SPDEwf}
\eea  
so we see that~\eqref{eq:integrated} is quite natural. In particular, $\fs(t)$ is
a local/instantaneous wave speed and appears as a coefficient of the 
advective (spatial first-order derivative) term; cf.~\eqref{eq:secondorder}. 
Of course, one can always consider an averaging procedure for the wave speed
\benn
\frac{1}{t-t_0}\int_{t_0}^t \E[\fs(r)]~\txtd r,
\eenn
which is expected to be intimately related with time averages in~\eqref{eq:timelimits}. 
However, the key problem is that the SPDE~\eqref{eq:SPDEwf} itself cannot be 
solved \emph{alone}. There is currently no equation determining the \emph{unknown stochastic 
processes $\fs=\fs(t)$ or $p=p(t)$}. Due to the time-dependence of $\fs$, the differential 
equation~\eqref{eq:SPDEwf} remains an SPDE and does not reduce to an SODE. Therefore, the 
classical deterministic ODE/PDE arguments to determine the wave speed cannot be generalized 
directly. 

Of course, re-centering can also be viewed more abstractly as an optimization problem, where we 
try to adapt our moving frame to minimize the distance to a some general reference solution. This
approach can be formalized via the optimization problem
\be
\label{eq:minLT}
\min_{y\in\R} \cM(y),\qquad \cM(y):=\|u(t)-u_{\textnormal{ref}}(\cdot-y)\|,
\ee
where two choices have to be made. Firstly, the function $u_{\textnormal{ref}}$ is the 
given \emph{reference solution} often still taken as the family 
of deterministic waves $u_{\textnormal{ref}}=\Phi_y=\Phi(\cdot -y)$; see~\cite{LordThuemmler}. 
Secondly, one has to select the type of (spatial) norm $\|\cdot\|$ to be
used in~\eqref{eq:minLT}. A common choice is just the $L^2(\R)$-norm~\cite{LordThuemmler} or
a weighted version $L^2(\R;\rho)$ with a weight function $\rho=\rho(x)$, which one could also adapt
to the position of the wave $p(t)$; see~\cite{KruegerStannat1}. In any case, let us agree
from now on that we have chosen a Hilbert space so that $\|\cdot\|=\langle\cdot,\cdot\rangle$,
where $\langle\cdot,\cdot\rangle$ denotes the inner product. The minimization then 
implicitly defines the position of the wave $p=p(t)$ by computing a minimizer $y$ for each 
$t\geq 0$ and setting $y=:p(t)$. However, one first has to guarantee that a minimizer 
exists~\cite{InglisMacLaurin}. 
Although it does exist under reasonable conditions, it does not have to be unique. It
is easy to imagine that there can be jumps for the \emph{global minimizer} so $p(t)$ may have 
jumps. Of course, we can always aim to track a \emph{local minimizer}. Differentiating $\cM(y)$ 
in $y$ yields a necessary \emph{critical point condition} and so we get the constraint 
\be
\label{eq:cpcond}
\langle \Phi_y', u(t)-\Phi_y \rangle=0.
\ee
If the second derivative of $\cM(y)$ 
\be
\frac{\txtd^2 }{\txtd y^2} \cM(y)=-2\langle u(t),\Phi''_y\rangle,
\ee
vanishes, the local minimizer becomes degenerate, so we have a criterion to test for 
potential jumps. Now, we have essentially given another implicit 
\emph{definition}~\cite{LordThuemmler,InglisMacLaurin,KruegerStannat1}
of the wave speed for a stochastic travelling wave and we have chosen $p=p(t)$ as moving-frame 
coordinates for its measurement. Next, one may ask, whether there is a differential
equation for $p(t)$? Starting from the form~\eqref{eq:minLT}, we may use~\eqref{eq:cpcond} 
and a variant of It\^o's formula 
to derive a stochastic ordinary differential equation (SODE) for the position up to the 
first jump time of $p(t)$. Consider the It\^o-SPDE~\eqref{eq:SPDEmath1} and suppose we have 
additive noise $g(u)\equiv\sigma$ for simplicity. Then a relatively direct 
calculation~\cite{InglisMacLaurin} yields
\bea
\label{eq:SDEwavespeed}
\txtd p &=& \left[\frac{\langle w,\partial_x^2 \Phi'_p\rangle+\langle f(w+\Phi_0)-
f(\Phi_0),\Phi'_p\rangle}{\zeta(p,w)} 
+ \sigma^2\frac{\langle Q\Phi'_p,\Phi''_p\rangle}{\zeta(p,w)^2} \right.\nonumber\\
&&+ \left.\frac{\sigma^2}{2}\frac{\langle w+\Phi_0-\Phi_p,\Phi'''_p\rangle \langle Q\Phi'_p,
\Phi'_p\rangle}{\zeta(p,w)^3} \right]~\txtd t+\frac{\sigma}{\zeta(p,w)}~\langle \Phi'_p,
\txtd W \rangle 
\eea
where $p=p(t)$, and
\benn
w=w(t):=u(t)-\Phi_0,\quad \zeta(p(t),w(t)):=\langle\Phi_0',\Phi_{p(t)}'\rangle-
\langle w(t),\Phi_{p(t)}''\rangle.
\eenn
It is crucial to note that~\eqref{eq:SDEwavespeed} is now an SODE for the wave speed. 
Unfortunately, the SODE is rather involved and usually impossible to solve analytically, 
even if an explicit representation of the deterministic wave profile is available. This difficulty 
is completely expected from the classical analysis of deterministic travelling wave problems (see 
Section~\ref{sec:PDE}) since we are essentially trying to capture the deviation of the stochastic
solution from a deterministic reference. Structurally, this amounts to studying an evolution 
equation for perturbations, and we know that the deterministic stability problem always involves
the shape of the travelling wave as an implicit input, even in the linear stability 
problem~\eqref{eq:linPDEtw}. The same effect precisely occurs in the position 
SODE~\eqref{eq:SDEwavespeed}. 

We emphasize that one also finds several variants of dynamic re-centering in the 
literature~\cite{CartwrightGottwald,Dierckxetal}, which can lead to different 
evolution equations, e.g., using different spatial norms and/or other adapted
coordinate systems. For example, one could consider~\eqref{eq:minLT}, fix 
a deterministic wave $\Phi$ and subtract a deterministic speed~\cite{KruegerStannat1},
and still consider a minimization approach via
\be
\label{eq:minLT1}
\min_{y\in\R} \|u(t)-\Phi(\cdot-st-y)\|=:\min_{y\in\R} \cM_s(y).
\ee
As another modification, we could minimize~\eqref{eq:minLT1} by \emph{relaxation} 
considering an evolution equation along the negative gradient of $\cM_s(y)$. In this 
approach we never solve the optimization problem directly but use relaxed gradient 
descent dynamics induced by the functional, which is to be 
minimized. More precisely, fixing a relaxation parameter $r>0$ one then has to 
consider the (pathwise-defined) random ordinary differential equation (RODE) 
\be
\label{eq:yrcoupled}
\frac{\txtd y^r}{\txtd t} = -r \langle u(t)-\Phi(\cdot-st-y^r),
\partial_x \Phi(\cdot-st-y^r)\rangle,\qquad \fs^r(t):=\frac{\txtd y^r}{\txtd t}
\ee
where $y^r=y^r(t)$ has an index $r$ to indicate its dependence upon \emph{imposed} 
relaxation parameter and $\fs^r$ measures the deviation of the wave speed since 
it is the derivative of the position deviation $y^r(t)=y^r$ from the deterministic 
profile. Now define $p^r(t):=-st-y^r(t)$ and let 
\benn
u^r(t):=u(t)-\Phi(\cdot + p^r),\qquad p^r=p^r(t),
\eenn
so $u^r=u^r(t)$ essentially measures the fluctuations of the SPDE around a dynamically 
phase-adapted deterministic wave. For additive noise $g(u)\equiv \sigma$, the 
evolution equation of this deviation is then given by~\cite{KruegerStannat1}
\be
\label{eq:KruegerStannat1}
\txtd u^r = \left[\partial_x^2 u^r + f(u^r+\Phi(\cdot+p^r))-f(\Phi(\cdot +p^r))+
 \frac{\txtd y^r}{\txtd t} \partial_x \Phi(\cdot + p^r)\right]
+\sigma~\txtd W.
\ee
Using It\^o's formula one can now also show that $\fs^r(t)$ satisfies an 
SODE~\cite[Lem.~3.2]{KruegerStannat1} since  
\benn
\fs^r(t)=\frac{\txtd y^r}{\txtd t}=-r \langle u^r(t),
\partial_x \Phi(\cdot-st-y^r)\rangle
\eenn
and we already have an equation for $u^r$ given by~\eqref{eq:KruegerStannat1}. Yet, the 
evolution equation for $\fs^r(t)$ is again fully coupled to the original SPDE. Therefore, 
even though we tried to simplify the procedure by not looking at the position $p$ and speed
$s$ but at the relaxed version $y^r$ of the position and deviation of the 
speed $\fs^r(t)$, we still face a (set of) stochastic nonlinear evolution equation(s). 
In view of this challenge, a simplification is to assume that one is 
interested in the case of \emph{small noise} deviations from a deterministic profile
\benn
0<\sigma\ll1,\qquad u(0)=\Phi_0, \qquad w=u-\Phi_0,    
\eenn
so that $w=w(t)$ is a stochastic process measuring the deviation from the deterministic
solution as for $\sigma=0$ we would have $u(t,x)=\Phi(x-st)=\Phi_0(\eta)$. We shall 
encounter \emph{multiscale} results based upon small-noise assumptions 
in Section~\ref{ssec:bistoch}.

%%%%%%%%%%%%%%%%%%%%%%%%%%%%%%%%%%%%%%%%%%%%%%%%%%%%%%%%%%%%%%%%%%%%%%%%%%%%%%%%%%%%%%%%%%%%%%
\subsection{Monostable Stochastic Waves}
\label{ssec:mostoch}

Here we collect results regarding different variants of the monostable FKPP-SPDE
\be
\label{eq:SPDEFKPP1}
\txtd u = \left[\partial_x^2 u + u(1-u)\right]~\txtd t + g(u)~\txtd W,\qquad u_0(x)=u(0,x).
\ee
The first question is whether we can prove that solutions resembling travelling waves 
can actually exist. Using a long-range voter model~\cite{MuellerTribe1} and/or the duality
to particle systems~\cite{DoeringMuellerSmereka}, one natural setting to consider is 
\be
\label{eq:MuellerSowers}
g(u)=\sigma \sqrt{u(1-u)},\qquad Q=\Id,\qquad \text{$\sigma>0$ sufficiently small.}
\ee
Suppose the initial condition is locally supported close to a front according 
to~\eqref{eq:intlocsupp}. It is then proven in~\cite{MuellerSowers} for~\eqref{eq:SPDEFKPP1}
with noise term~\eqref{eq:MuellerSowers} that defining $a(t),b(t)$
according to~\eqref{eq:abdef}, then we almost surely have that the limit  
\benn
\lim_{t\ra +\I} \frac{b(t)}{t}=:b_*(\sigma)
\eenn
exists and $b_*(\sigma)$ is non-random; see also~\cite{ConlonDoering,OksendalVageZhao,
OksendalVageZhao1}. In addition, the law of $u(t,b(t)+x)$ tends to a stationary limit as 
$t\ra +\I$ as does the length distribution of the interval $[a(t),b(t)]$. This means we 
really obtain a true 
front-type wave as $t\ra +\I$. Then one may ask, how $b_*$ depends upon $\sigma$.
Based on formal physical approximation arguments and numerical evidence, it has been 
found~\cite{BrunetDerrida,BrunetDerrida2} (see also \cite{BreuerHuberPetruccione,
KesslerNerSander,PechenikLevine}) that
\be
\label{eq:BrunetDerrida}
b_*(\sigma)=2-\frac{\pi^2}{|\ln\sigma^2|^2}+\text{higher-order terms}\qquad \text{as 
$\sigma\ra 0^+$},
\ee
which is also known as the \emph{Brunet-Derrida conjecture}. Note that the wave moves slower than
the classical deterministic Fisher wave. The inverse-logarithmic correction 
in~\eqref{eq:BrunetDerrida} is much larger than expected by a naive asymptotic
expansion. Of course, it just results from the pushed nature of the front, i.e., from the 
interplay between small noise and instability~\cite{KuehnCurse} near the unstable leading 
edge of the front. More precisely, if we are sufficiently close to the leading edge, say 
$u<\sigma^2$, then $u(1-u)<\sigma\sqrt{u(1-u)}$, so the noise term dominates the reaction 
term. The Brunet-Derrida conjecture~\cite{BrunetDerrida2} has been formed by relating the 
speed to the FKPP with a \emph{cut-off}
\be
\label{eq:cutoffmodel}
\partial_t v = \partial_x^2 v + v(1-v)\textbf{1}_{\{v\geq \sigma^2\}},
\ee
which can be studied using elegant deterministic 
arguments~\cite{DumortierPopovicKaper,DumortierPopovicKaper1} based upon geometric 
singular perturbation theory~\cite{KuehnBook}. The front for the cut-off
system does behave rather ``weakly pushed''~\cite{PanjavanSaarloos}. However, even with the 
results for~\eqref{eq:cutoffmodel}, one still has to rigorously connect 
the cut-off model~\eqref{eq:cutoffmodel} with the SPDE~\eqref{eq:SPDEFKPP1} or with an 
underlying particle system. The full Brunet-Derrida conjecture has been proven 
in~\cite{MuellerMytnikQuastel} (see 
also~\cite{BerardGouere}), i.e., the small-noise expansion of the stochastic travelling 
front wave speed for~\eqref{eq:SPDEFKPP1}-\eqref{eq:MuellerSowers}
is indeed given by~\eqref{eq:BrunetDerrida} with higher-order terms 
of order $\cO((\ln|\ln\sigma|)|\ln\sigma|^{-3})$; see also~\cite{ConlonDoering} for a 
lower bound. The same wave speed asymptotics also hold rigorously for
the noise term $g(u)=\sigma\sqrt{u}$, which can be derived from a long-range contact 
process~\cite{MuellerTribe1}. Indeed, it has been proven~\cite{Tribe} that travelling 
wave-type solutions also exist for
\be
\label{eq:Tribesimsq}
g(u)=\sqrt{u},\qquad u_0(x)=\min(1,\max(-x,0)),\qquad Q=\Id,\qquad \text{$\sigma>0$ 
sufficiently small.}
\ee
with a well-defined limit $\lim_{t\ra +\I} b(t)/t$; note that $a(t) =-\I$ due 
to the form of the noise term effectively acting like additive noise near $u=1$. 
In addition to looking at speeds, 
one can also define a diffusion coefficient of the front~\cite{BrunetDerrida2}
\benn
D(\sigma):=\lim_{t\ra +\I}\frac{\E[b(t)^2]-\E[b(t)]^2}{t},
\eenn 
yet its asymptotics for monostable equations seems to be more difficult to analyze 
rigorously~\cite{BrunetDerridaMuellerMunier,Panja1,RoccoCasademuntEbertvanSaarloos}; 
we remark that the diffusion coefficient seems to have an interesting relation to 
coalesence times in the context of related particle models~\cite{BrunetDerridaMuellerMunier1}. 

\begin{figure}[htbp]
	\centering
	\begin{overpic}[width=1\textwidth]{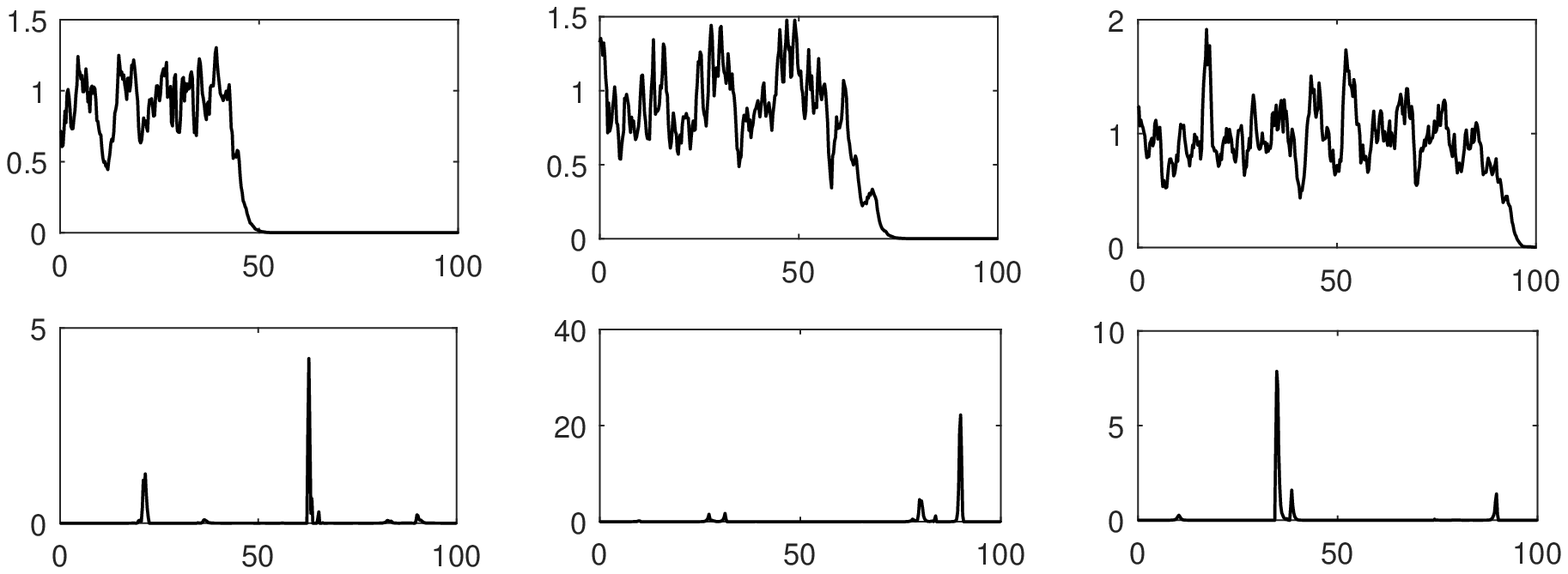}
		\put(0,22){\scalebox{1}{$\underline{\sigma=1}$}}
		\put(0,6){\scalebox{1}{$\underline{\sigma=3}$}}
		\put(16,0){\scalebox{0.7}{$x$}}
		\put(45,0){\scalebox{0.7}{$x$}}
		\put(73,0){\scalebox{0.7}{$x$}}
		\put(16,15){\scalebox{0.7}{$x$}}
		\put(45,15){\scalebox{0.7}{$x$}}
		\put(73,15){\scalebox{0.7}{$x$}}
		\put(9,10){\scalebox{0.7}{$u$}}
		\put(9,26){\scalebox{0.7}{$u$}}
		\put(37,10){\scalebox{0.7}{$u$}}
		\put(37,26){\scalebox{0.7}{$u$}}
		\put(65,10){\scalebox{0.7}{$u$}}
		\put(65,26){\scalebox{0.7}{$u$}}
	\end{overpic}
	\caption{\label{fig:05}Numerical simulation of~\eqref{eq:propfail} with white noise
	truncated at 25 eigenmodes as above on a domain $[0,100]$ with the same initial 
	condition for both rows. The top row shows $\sigma=1$, 
	which still exhibits a stochastic travelling wave, while propagation failure occurs in 
	the bottom row for $\sigma=3$. The solution is displayed at the same time snapshots each 
	column.}
\end{figure}

The next natural question is, what 
happens if we consider large(r) noise. Formal approximations and numerical simulations 
suggest that the wave speeds change substantially~\cite{HallatschekKorolev}. Even more
drastically, if the noise is large enough, we may have \emph{propagation failure}. For 
example, consider the suitably scaled monostable SPDE
\be
\label{eq:propfail}
\txtd u = \left[\partial_x^2 u +\sigma^2u(1-u)\right]~\txtd t + 
\sigma^2 u~\txtd W,\qquad Q=\Id,~u_0(x)\geq 0. 
\ee
Then it is proven~\cite{ElworthyZhaoGaines} that for small noise a wave exists. However, 
given a suitably fixed initial condition, we can pick any compact set $\cK\subset [0,\I)\times \R$ 
and there exists $\sigma_0>0$ sufficiently large that for any $\sigma\in[\sigma_0,+\I)$ we have
\be
\label{eq:fail1}
\P\left(\sup_{(t,x)\in\cK}u(t,x)>\txte^{-K_1t\sigma^4}\right)<\txte^{-K_2t^2\sigma_0^4}
\ee  
for some constants $K_{1,2}>0$. So the solution is exponentially small with 
high probability~\cite{ElworthyZhaoGaines}; see also~\cite{DaviesTrumanZhao,ElworthyZhao};
see also Figure~\ref{fig:05}.

In theoretical terms this implies $u(t,x)\ra 0$ as $\sigma\ra +\I$ and/or $t\ra +\I$, while 
in practical terms it means that we can expect that a travelling wave-form initial condition 
is going to lose its wave-form shape. This propagation failure is observed in simulations for 
initial data with compact support such as approximate identities (``approximate $\delta$-distributions'') 
and strong noise~\cite{ElworthyZhaoGaines,Gaines,KuehnFKPP}, while in the deterministic case 
$\sigma= 0$ we know there is one left-moving and one right-moving front for such 
initial conditions. The dichotomy of propagation failure was made even more precise~\cite{MuellerTribe2} 
for the SPDE
\benn
\txtd u = \left[\partial_x^2 u + \theta u-u^2 \right] +\sqrt{u}~\txtd W,\qquad 
Q=\Id,
\eenn
for suitably chosen non-negative compactly supported initial data. In this case, one can 
prove~\cite{MuellerTribe2} that there exists a constant $\theta_\txtc$ independent of $u_0$ 
such that
\be
\label{eq:fail2}
\P(u(t,0)\neq 0~\forall t>0)~\left\{
\begin{array}{ll}
=0&\text{ if $\theta<\theta_\txtc$,}\\
>0&\text{ if $\theta>\theta_\txtc$.}
\end{array}
\right.
\ee
So the process dies out if the local linear deterministic instability induced by the
reaction term is not strong enough. In fact, the existence of travelling waves in the 
survival regime holds, i.e., we may replace the small-noise condition 
in~\eqref{eq:Tribesimsq} by $\theta>\theta_\txtc$. We refer 
to~\cite{Tribe,Kliem,Kliem1,HorridgeTribe} for more details on limiting distributions 
and the role of initial conditions in this case, while the effect of large noise
is also discussed in~\cite{MuellerMytnikRyzhik}.  
Of course, propagation failure effects such as~\eqref{eq:fail1} 
and \eqref{eq:fail2} can also occur for various other noises~\cite{HuangLiu}.

We remark that there is strong numerical evidence that the statistical properties of 
a stochastic wave even contain early-warning signs to indicate closeness to propagation 
failure~\cite{KuehnFKPP}; see also Section~\ref{sec:outlook}. 
In addition, there is growing evidence that suitably chosen noise can arrest/freeze pulled 
fronts or change their direction~\cite{MeersonSasorov}. For viewpoints relating stochastic 
monostable dynamics to finite-size effects in autocatalytic reactions we refer 
to~\cite{LemarchandLesneMareschal,MaiSokolovBlumen,Kuzovkovetal}. Another possible direction is to 
consider local modifications of the monostable reaction term to balance certain noise-induced 
effects~\cite{PanjavanSaarloos1} or to study the transition between pulled and pushed noisy
waves~\cite{BirzuHallatschekKorolev,RoccoRamirez-PiscinaCasademunt}. Now we are also going to 
transition and change from the monostable to the bistable case.

%%%%%%%%%%%%%%%%%%%%%%%%%%%%%%%%%%%%%%%%%%%%%%%%%%%%%%%%%%%%%%%%%%%%%%%%%%%%%%%%%%%%%%%%%%%%%%
\subsection{Bistable Stochastic Waves}
\label{ssec:bistoch}

\begin{figure}[htbp]
	\centering
	\begin{overpic}[width=1\textwidth]{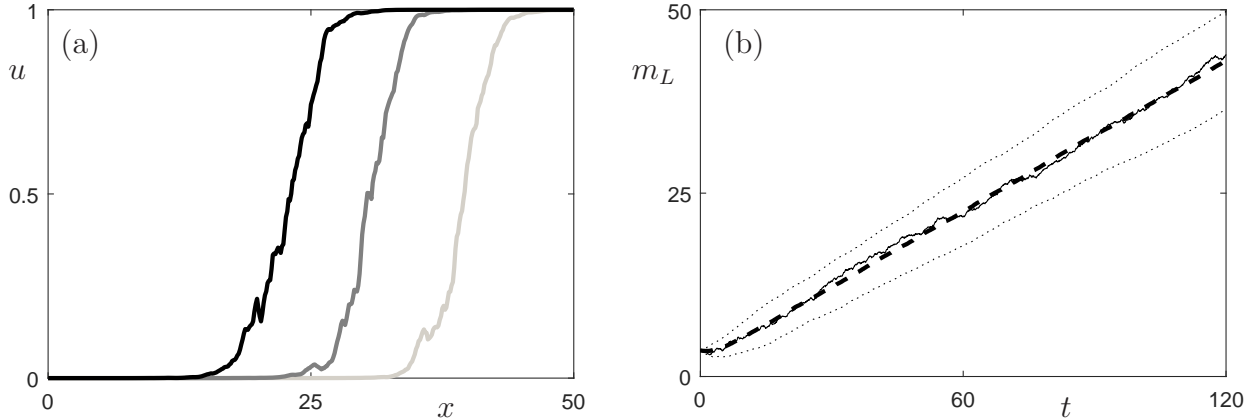}
		\put(6,28){\scalebox{1}{(a)}}
		\put(57,28){\scalebox{1}{(b)}}
		\put(83,0){\scalebox{1}{$t$}}
		\put(35,0){\scalebox{1}{$x$}}
		\put(50,26){\scalebox{1}{$m_L$}}
		\put(2,26){\scalebox{1}{$u$}}
	\end{overpic}
	\caption{\label{fig:06}(a) Direct numerical simulation of the bistable 
	case~\eqref{eq:SPDEbi1} with multiplicative noise~\eqref{eq:noise1Nag}. The 
	parameters are $a=\frac14$ and $\sigma=1.5$. The spatial domain is $[0,50]$
	and the same white noise truncated after 25 modes is used as above for the 
	monostable simulations. (a) A left-moving stochastic wave is shown at three
	time snapshots. (b) The corresponding $\tilde{m}_L(t)$ for $L=50$ is shown as
	a solid curve, while the average $\E[\tilde{m}_L(t)]$ is the dashed thicker
	black curve. The dots indicate two standard deviations from the mean, i.e., the 
	two dotted curves are $\E[\tilde{m}_L(t)]\pm 2\sqrt{\text{Variance}(\tilde{m}_L(t))}$,
	where mean and variance have been computed over $100$ sample paths.}
\end{figure}

It has been known for a quite a long time based on physical grounds and simulations 
that in the bistable case $f=f_3$, front-like solutions exist, at least in the weak noise 
setting~\cite{MikhailovSchimanskyGeierEbeling}. For these fronts one may ask similar 
questions regarding the wave speed of propagating fronts as for the monostable case. 
In~\cite{Armeroetal} a noise term is derived from external fluctuations in the 
control parameter $a$ to yield the model
\be
\label{eq:SPDEbi1}
\txtd u = \partial_{x}^2 u + u(1-u)(u-a)~\txtd t+\sigma \tilde{g}(u)\circ\txtd W,
\ee
where $W$ is white in time, has spatial correlation function $q=q(x-y)$, and 
the noise term is to be interpreted in the Stratonovich form~\cite{Panja} in this
context. Examples for $\tilde{g}$ will be discussed below. A typical 
\emph{formal approximation} approach~\cite{Armeroetal,Armeroetal1} 
for~\eqref{eq:SPDEbi1} is to reduce the problem to moments. One often finds a 
reference~\cite{Armeroetal,Armeroetal1} to ``Novikov's Theorem''~\cite{Novikov} in 
the theoretical physics literature, which is used to derive approximating equations; 
see also~\cite{Panja,GarciaOjalvoSancho,SantosSancho} for some further references using 
this approach. Interestingly, the original reference by Novikov~\cite{Novikov}
frequently cited in this context of the stochastic travelling waves literature 
in the physical sciences does not contain any theorems but only proposes a
formal approximation based on physical grounds for noise in the context of turbulence
problems. The name ``Novikov'' and the transformation to a different stochastic
process appearing in the theorem, may lead one intuitively to think of
the so-called Novikov condition known to appear in he context of Girsanov's 
Theorem~\cite{Protter}. Yet, the the Novikov condition is by A.A.~Novikov~\cite{Novikov1}
while ``Novikov's Theorem'' is by E.A.~Novikov~\cite{Novikov}. In fact, it is probably
better to use the convention Furutsu-Novikov Theorem~\cite{Pecseli,KonotopVazquez} for
the latter to avoid confusion\footnote{I would like to thank Eulalia Nualart for pointing 
out the alternative attribution to Furutsu~\cite{Furutsu} to me, which I had not 
been previously aware of.}. The Furutsu-Novikov Theorem states that
\be
\label{eq:NovikovFurutsu}
\E\left[\Xi(t)\cF[\Xi\left(\left.\cdot\right|_{t_0}^t\right)]\right]=
\int_{t_0}^t C_{\Xi(\cdot)\Xi(\cdot)}(t,s)\E\left[\frac{\delta 
\cF\left[\Xi(\cdot|_{t_0}^t)\right]}{
\delta \Xi(s)}\right]~\txtd s,
\ee
where $\Xi$ is any zero-mean Gaussian process with given \emph{autocovariance 
function} $C_{\Xi\Xi}(t,s)$, $\cF\left[\xi(\cdot|_{t_0}^t)\right]$ is a functional of $\xi$ 
over the time interval $[t_0,t]$ and $\delta/\delta\xi(s)$ denotes the functional 
derivative with respect to $\xi$ at $s$. One idea to utilize this theorem\footnote{The 
terminology of 'theorem' is potentially not ideal. Although the formal relation seems 
evident from a calculation, most justifications in the literature of the Furutsu-Novikov Theorem 
seem to be extremely concise in terms of their explanation and formalization.} for travelling 
waves is to notice that direct averaging of~\eqref{eq:SPDEbi1} will not produce 
a deterministic PDE to leading-order as $\sigma \tilde{g}(u)\circ \partial_tW$ 
may not have zero mean. Hence, one tries to restore this zero mean property. 
The mean value of the noise can be calculated via the 
Furutsu-Novikov Theorem~\cite{Armeroetal,Armeroetal1}
\be
\label{eq:ItoNov}
\sigma \E[\tilde{g}(u)~\circ \partial_t W]=\frac{\sigma^2}{2}
q(0)\E[\tilde{g}'(u)\tilde{g}(u)],
\ee
where we observe that the right-hand side of formula~\eqref{eq:ItoNov} is just the 
average of an It\^o-Stratonovich correction term. Re-writing the 
SPDE~\eqref{eq:SPDEbi1} suggestively with $\partial_t W=:\xi$ as
\bea
\partial_t u &=& \partial_x^2 u + f_3(u) + \frac{\sigma^2}{2} q(0)\tilde{g}'(u)\tilde{g}(u)+ 
\tilde{g}(u)\xi - \frac{\sigma^2}{2} q(0)\tilde{g}'(u)\tilde{g}(u),\nonumber \\
&=:& \partial_x^2 u + f_3(u) + \frac{\sigma^2}{2} q(0)\tilde{g}'(u)\tilde{g}(u) + R,\qquad R=R(t,x,u) 
\label{eq:funnyphysics}
\eea
means that if we average now, the new noise term $R$ should disappear due to zero
mean. Furthermore, the noise does have a changed correlation function
\benn
\E[R(t,x,u)R(s,y,u)]=\E[u(t,x)u(s,y)\xi(t,x)\xi(s,y)]+\cO(\sigma),\qquad \text{as 
$\sigma\searrow 0$}.
\eenn
Let us now illustrate, how to use~\eqref{eq:funnyphysics} for the case
\be
\label{eq:noise1Nag}
\tilde{g}(u)=u(1-u),
\ee
to formally track suitable averages for front-like structures. Fix a sufficiently 
big interval $[-L,L]$ within we suspect the front-like structure connecting $0$ to $1$. 
Next, define
\benn
\tilde{m}_L(t):=\int_{-L}^L u(t,x)~\txtd x
\eenn  
which is, yet another, useful random-variable to keep track of the position of the 
wave; see also Figure~\ref{fig:06}. 

Then define
the deviation of the position from the average as 
\benn
m_L(t):=\tilde{m}_L(t)-\E[\tilde{m}_L(t)]. 
\eenn
Note that this construction is yet another version of dynamic re-centering 
described in the first part of Section~\ref{sec:waves}. Setting now 
$u_{m_L}(t,x):=u(t,x+m_L(t))$ in~\eqref{eq:SPDEbi1} one obtains
\be
\label{eq:applyNovhere}
\txtd u_{m_L}=\left[\partial_x^2 u_{m_L} + \frac{\txtd m_L}{\txtd t} \partial_x u_{m_L} 
+ u_{m_L}(1-u_{m_L})(u_{m_L}-a)\txtd t \right] + u_{m_L}(1-u_{m_L})\circ \txtd W.
\ee 
Then one can evidently average again by defining the average front shape 
$u_{m_L}^0:=\E[u_{m_L}]$. What is the evolution equation for $u_{m_L}^0$? Taking
the average in~\eqref{eq:applyNovhere} and using the Furutsu-Novikov Theorem gives an
evolution equation for $u_{m_L}^0$. Yet, the first moment is coupled to higher
moments in general, leading to the problem of \emph{moment closure}~\cite{KuehnMC}. 
However, let us decompose the dynamics into two parts
\benn
u_{m_L}(t,x)= u_{m_L}^0(t,x) + u_{m_L}^\txtc(t,x),   
\eenn
where $u_{m_L}^0$ should describe the leading-order/mean dynamics of the re-centered front.
If one \emph{postulates}, say based upon numerical evidence, that the dynamics
of the deviations $u_{m_L}^\txtc$ decays quickly for any initial condition so
that $|u_{m_L}^\txtc|\ll 1$, then \emph{formally} keeping only lowest-order terms 
yields~\cite{Armeroetal,Armeroetal1}
\be
\label{eq:Armero}
\partial_t u^0_{m_L} = \partial_x^2 u^0_{m_L} + u^0_{m_L}(1-u^0_{m_L})
(c_0u^0_{m_L}-a_0),\qquad c_0=1-\frac{\sigma^2}{2}q(0),a_0=a+\frac{\sigma^2}{4} q(0). 
\ee
The equation~\eqref{eq:Armero} can now be analyzed using PDE techniques. For example,
we find that the wave speed changes depending upon varying $q(0)$ and $\sigma$. Yet,
the form of~\eqref{eq:Armero} already indicates that the wave speed change is far \emph{more
regular} for small noise compared to the monostable case, i.e., one observes a regular 
power-type expansion in $\sigma$, which is in accordance with results for bistable 
PDEs with a cut-off~\cite{BenguriaDepassierHaikala,KesslerNerSander,MendezCamposZemskov};
cf.~the discussion above regarding~\eqref{eq:cutoffmodel}.\medskip  

The preceding formal approximation~\eqref{eq:Armero} did not a-priori rely on 
direct \emph{small noise expansions} (or \emph{weak noise 
expansions})~\cite{DePasqualeGoreckiPopielawski}, which are another common tool, 
particularly for bistable equations~\cite{KaraziLemarchandMareschal,
MikhailovSchimanskyGeierEbeling,SchimanskyGeierMikhailovEbeling,SchimanskyZuelicke}. 
Of course, to analyze~\eqref{eq:Armero} as a meaningful approximation for stochastic
waves, then one still needs a small noise condition. In the small-noise regime, one
expects that if a travelling wave is deterministically locally stable, then the 
corresponding SPDE generates a similar wave-like profile with additional 
diffusive motion along the neutral (Goldstone) mode associated to the zero eigenvalue
generated by translation invariance; see also Section~\ref{sec:PDE}. These 
results can been made rigorous for bistable systems and fronts under certain 
assumptions~\cite{KruegerStannat1}. Suppose we start with the 
\be
\label{eq:SPDEbi2}
\txtd u = \partial_{x}^2 u + u(1-u)(u-a)~\txtd t+\sigma~\txtd W, \quad 
\textnormal{Tr}(Q)<+\I.
\ee
Then consider the SODE~\eqref{eq:yrcoupled} for the deviation of the position
$y^r(t)$ from the deterministic front and the SPDE~\eqref{eq:KruegerStannat1}
describing the fluctuations around the deterministic front $\Phi$; recall that both
quantities are dynamically re-centered. We want to study their equations on a finite
time interval, say $t\in[0,T]$. Assuming that $\sigma>0$
is sufficiently small, one may rigorously prove there is a decomposition
\benn
p^r(t)=-st-\sigma y_1^r(t)+o(\sigma),\quad \text{as $\sigma\ra 0$},
\eenn
where $s$ is the deterministic wave speed and 
\benn
y_1^r(t)=\int_0^t \fs^r_1(\tilde{t})~\txtd \tilde{t} 
\eenn
is the leading-order approximation of the position deviation. The 
leading-order approximation for the speed correction turns out to satisfy the 
SODE~\eqref{eq:KruegerStannat1}
\be
\label{eq:SODEloc}
\txtd \fs^r_1 = -r \fs^r_1 ~\txtd t - r\langle\partial_\eta 
\Phi(\cdot-st),\txtd W\rangle,
\ee 
where $\langle\cdot,\cdot\rangle=
\langle\cdot,\cdot\rangle_{L^2(\R,\rho(\cdot-st-\sigma y_1^r(t))}$ 
is a weighted norm moving within the adapted reference frame and one considers the
function $\rho(y)=Z\exp(s y)$, where $Z$ is constant computable from the
deterministic wave profile $\Phi$ as the normalization constant guaranteeing that
\benn 
\langle\txte^{c\cdot}\partial_\eta \Phi,\partial_\eta \Phi\rangle_{L^2(\R)}=1.
\eenn
Indeed, one may check that the function $\txte^{c\eta}\partial_\eta \Phi$ is the 
eigenfunction to the neutral eigenvalue zero of the adjoint $\cL^*$ to the 
operator $\cL$ arising from the linearization around the deterministic travelling wave 
defined in~\eqref{eq:linPDEtw}. Hence, we expect that in the weighted norm we just have
that the wave speed correction $\fs^r_1$ arises from a projection argument onto 
the eigenspace spanned by the neutral (Goldstone) mode spanned by $\partial_\eta \Phi$. 
This is indeed visible in the SODE~\eqref{eq:SODEloc} as there is deterministic decay 
from the relaxation parameter and diffusive wandering projected onto the neutral
mode as expected from physical intuition~\cite{Panja}. Furthermore, one may prove
a leading order SPDE approximation for the fluctuations $u^r(t)$ 
\benn
u^r(t)=\sigma u_1^r(t)+o(\sigma) ,\quad \text{as $\sigma\ra 0$},
\eenn
where the SPDE for $u_1^r(t)=u_1^r(t,x)$ is~\cite{KruegerStannat1}
\be
\label{eq:SPDEKS}
\txtd u_1^r=\left[\partial_x^2 u_1^r + f_3'(\Phi(\cdot -st))u_1^r
+\fs_1^r\partial_\eta\Phi(\cdot-st) \right]~\txtd t + \txtd W,\quad x\in\R.
\ee
The evolution equations~\eqref{eq:SODEloc} and~\eqref{eq:SPDEKS} only hold up to a 
stopping time $\tau=\tau(\sigma)$ as discussed already at the beginning of 
Section~\ref{sec:waves} since one may encounter jumps of the processes involved. Yet, 
we can always ensure for a given fixed maximal time $T>0$ that 
\benn
\lim_{\sigma\ra 0}\tau(\sigma)=T,
\eenn
so we truly have a small-noise approximation in a precise sense on finite 
time intervals. Although the equations~\eqref{eq:SODEloc} and~\eqref{eq:SPDEKS} are 
just the leading-order small-noise approximations one could in principle continue
the expansion in $\sigma$ yielding more evolution equations for higher-order 
corrections.

\begin{figure}[htbp]
	\centering
	\begin{overpic}[width=0.85\textwidth]{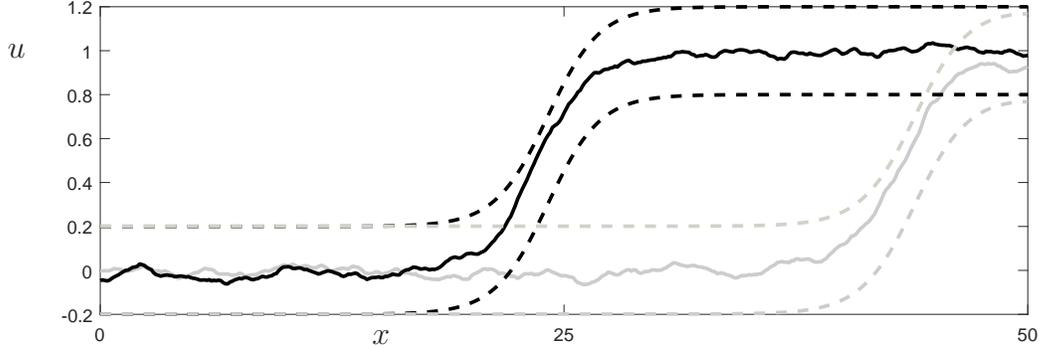}
		\put(35,0){\scalebox{1}{$x$}}
		\put(2,26){\scalebox{1}{$u$}}
	\end{overpic}
	\caption{\label{fig:07}Direct numerical simulation of the bistable 
	case~\eqref{eq:SPDEbi1} with additive noise and $\sigma=0.1$. The 
	parameters are $a=\frac14$ and $\sigma=1.5$. The spatial domain is $[0,50]$
	and the same white noise truncated after $25$ modes is used as above for the 
	simulations. A left-moving stochastic wave (solid curves) is shown at two
	time snapshots. Furthermore, we show two neighbourhoods of the wave (dashed
	curves) constructed from the deterministic profile. These neighbourhoods
	are useful in bounds such as~\eqref{eq:Sprob}.}
\end{figure}

The next natural question one can pose is, how can we analyze approximating equations,
either without a small-noise or with a small-noise assumption? Since these evolution
equations are often simpler, yet still impossible to solve explicitly, one aims for 
estimates~\cite{Stannat}. For example, consider the SPDEs for the diffusion
along the neutral mode~\eqref{eq:KruegerStannat1} or its leading-order 
approximation~\eqref{eq:SPDEKS} and the stopping times
\benn
\tau_{K}:= \inf \{t\geq 0:\|u^r(t)\|>K\}\quad \text{and}\quad 
\tau^1_{K}:= \inf \{t\geq 0:\|u^r_1(t)\|>K\},
\eenn
for some fixed spatial norm $\|\cdot\|$. These stopping times provide qualitative 
information about the times when fluctuations become larger
than a given constant $K>0$. In particular, one key aim is to prove bounds on
the distributions 
\be
\label{eq:Sprob}
\P(\tau_K\leq \kappa)\quad \text{ and }\quad\P(\tau^1_K\leq \kappa). 
\ee
Alternatively, we could also look directly at the probabilities 
\be
\label{eq:HHprob}
\P\left(\sup_{t\in[0,T]}\|u^r(t)\|>\kappa\right)\quad \text{and}\quad 
\P\left(\sup_{t\in[0,T]}\|u^r_1(t)\|>\kappa\right),
\ee
for some $\kappa>0$; see also Figure~\ref{fig:07}. 
Of course, similar remarks apply to probability bounds on stochastic 
corrections to the deterministic speed or position of the wave. 
It is then quite natural to expect that these probabilities can
be estimated using arguments such as Markov's inequality, Chebyshev's
inequality, Doob's inequality, or the Burkholder-Davis-Gundy and other 
related concentration inequalities for probabilities~\cite{Kallenberg}. These
inequalities often provide an elegant way to convert the question to deterministic 
arguments about moment bounds. We may view estimates on the probabilities
\eqref{eq:Sprob} and/or \eqref{eq:HHprob} of the fluctuations near a deterministic wave
or related results on the speed/position as \emph{stability results}, which have been 
studied already in quite some detail in the bistable 
setting~\cite{HamsterHupkes,InglisMacLaurin,KruegerStannat1,Stannat,Stannat1}. 

To see that we indeed expect stability for small noise, we state at least one
result in this direction~\cite{HamsterHupkes} in a very special case. Let $\Phi$ be 
the deterministic front for $f=f_3$ with speed $s$. Consider the SPDE
\be
\label{eq:HH}
\txtd u = \left[\partial_x^2 u +u(1-u)(u-a)\right]~\txtd u+\sigma\tilde{g}(u)~\txtd B,\quad
u(0,x)=u_0(x),
\ee
where $B=B(t)$ is a standard real-valued Brownian motion and
$\tilde{g}(\Phi)=-\sqrt2 \partial_x \Phi$. Hence, the noise only acts ``rigidly'' on the 
stochastic wave as we have enforced an invariance along the neutral (Goldstone) mode 
for $g$. Upon starting with a well-prepared initial condition, which is computable just
using $\Phi(x)$, one
can actually ensure that we only see a stochastic wave with changed speed and fluctuating 
position but no additional fluctuations around the wave, so fluctuation estimates 
are trivial. Furthermore, (yet another) variant of the definition of position, 
say $p_*(t)$, is shown to satisfy~\cite{HamsterHupkes}
\be
\label{eq:nicebi}
p_*(t)=\frac{st}{\sqrt{1+\sigma^2}}+\frac{\sigma}{\sqrt{1+\sigma^2}}B(t).
\ee
Hence, this explicit (very special-case!) formula now easily yields stability type results
as we only need to estimate the position and speed, which depend in a simple way 
on Brownian motion so well understood upper/lower bounds for Brownian motion 
can be applied~\cite{MoertersPeres}. It is far more difficult to obtain general 
stability results but the bistable case, as illustrated by formula~\eqref{eq:nicebi}, 
is expected to be quite tame
in the small noise regime~\cite{HamsterHupkes,HamsterHupkes2,InglisMacLaurin,
KruegerStannat1,Stannat,Stannat1}; cf.~formula~\eqref{eq:BrunetDerrida} for the 
monostable case. 

Lastly, we point out that one should always keep in mind that also for the 
bistable case, any theoretical result should be compared to microscopic modelling 
of the noise~\cite{MeersonSasorovKaplan,KhainMeerson}; see also~\cite{Panja}.

%%%%%%%%%%%%%%%%%%%%%%%%%%%%%%%%%%%%%%%%%%%%%%%%%%%%%%%%%%%%%%%%%%%%%%%%%%%%%%%%%%%%%%%%%%%%%%
\section{Summary \& Outlook}
\label{sec:outlook}

There are many topics closely connected to travelling waves for monostable and bistable 
SPDEs. We mention a few of these directions here. In fact, there is an even simpler SPDE,
which can generate interface-like solutions~\cite{Tribe1} given by
\be
\label{eq:Tribesimple}
\txtd u = \partial_x^2 u~\txtd t + \sigma\sqrt{u(1-u)}~\txtd W,\qquad Q=\Id.
\ee
The interfaces of~\eqref{eq:Tribesimple} behave like Brownian motion. The 
model~\eqref{eq:Tribesimple} can be derived from a long-range voter 
model~\cite{MuellerTribe1,Kliem2} and is therefore microscopically related to 
the monostable FKPP equation~\eqref{eq:SPDEFKPP1}-\eqref{eq:MuellerSowers}, which 
has the same noise term. 

Another topic related to the FKPP equation is its generalization to higher dimensions 
\be
\label{eq:FKPPhigh}
\txtd u = \left[\Delta u + u(1-u)\right]~\txtd t + \sigma g(u)~\txtd W,\qquad u=u(t,x),~x\in\R^d,
\ee
for some $d\geq 2$. Upon using a certain multiplicative noise and a suitable initial 
condition, one can again obtain propagating front-like solutions invading the deterministically
unstable state~\cite{Moro,RiordanDoeringBenAvraham}. Since the interface propagation has
now a non-trivial spatial structure~\cite{RoccoRamirez-PiscinaCasademunt,TripathyvanSaarloos,
Tripathyetal}, e.g., a curve-like rough interface for $d=2$, it is
natural to try to connect its dynamics to effective interface models such as mean 
curvature flow, the Mullins-Skekerka equation, or the Kardar-Parisi-Zhang (KPZ) equation.
The KPZ equation is given by
\be
\label{eq:KPZ}
\partial_t h = \Delta h + (\nabla h)^2+\partial_t W,\qquad Q=\Id,~h=h(t,x),~x\in\R^{d-1}.
\ee 
The KPZ equation~\cite{KardarParisiZhang} is one normal-form type model or 
``universality class'' for interface growth, where one can think of $h$ as a height 
function of the interface for $d=2$. However, note that~\eqref{eq:KPZ} is not well-posed
as written in the form~\eqref{eq:KPZ} since the regularity of the space-time 
white-noise $\xi=\partial_t W$ does not allow one to define $(\nabla h)^2$ via a standard
fixed-point argument to obtain the existence of solutions. Due to this regularity issue, 
the KPZ equation is an example of a \emph{singular SPDE}. Under certain technical assumptions, 
it can be possible to \emph{renormalize} a singular SPDE and analyze it within 
the framework of regularity structures~\cite{Hairer1,Hairer2} or within paracontrolled 
distribution theory~\cite{GubinelliImkellerPerkowski,GubinelliPerkowski}. Of course,
if we view the monostable and bistable SPDEs in higher dimensions and/or with very 
irregular noise terms, they need renormalization as well, which has been noticed 
already in the context of numerical simulation; see~\cite{Panja} and references therein. 

Instead of considering higher spatial dimensions $d\geq 2$ for scalar equations, one
may ask, what happens if $d=1$ but we consider systems of reaction-diffusion SPDEs with
various nonlinearities? The theory for travelling waves in this context is even less 
developed. The typical results/effects for the one-component case are still key points 
for systems, e.g., front-like structures including propagation 
failure~\cite{MuellerTribe} or stability results~\cite{HamsterHupkes1} have been proven for
two-component model problems. However, there are additional new phenomena possible if we 
consider systems for $d\geq 2$ such as 
spiral-like structures~\cite{Sendinaetal,Dierckxetal}. It is natural to conjecture that 
spiral-like waves can be found if we perturb the classical models for spiral waves
such as the FitzHugh-Nagumo~\cite{FitzHugh,Nagumo} equation, the Barkley 
model~\cite{Barkley2}, or the Oregonator system by 
noise~\cite{BerglundKuehn,BonaccorsiMastrogiacomo,LaingLord,SauerStannat,Sendinaetal}.

Instead of generalizing to higher spatial dimensions, there are also first 
attempts to consider waves for other noise terms, such as L\'evy noise~\cite{BrockmannHufnagel}.
Furthermore, one may replace the heat equation part $\partial_t u=\Delta u$ by 
more general fractional derivative operators~\cite{BrockmannHufnagel} derived from
anomalous diffusion, or even observe anomalous diffusion from classical 
equations~\cite{RoccoEbertvanSaarloos}. In fact, the analysis of travelling waves
for deterministic PDE involving fractional operators is another recently emerging
area for monostable~\cite{CabreRoquejoffre,del-Castillo-NegreteCarrerasLynch,Engler} 
as well as bistable~\cite{AchleitnerKuehn1,AchleitnerKuehn2,Chmaj,GuiZhao,
VolpertNecNepomnyashchy,Zanette} cases. One should even suspect that 
waves for fractional diffusion operators for PDEs and waves for 
SPDEs are deeply connected~\cite{RoccoEbertvanSaarloos} since both underlying
classes of differential equations are derived from very similar microscopic stochastic 
modelling principles.

Nonlocal fractional operators are just one class, where nonlocality has recently entered
into focus. Another important recent class motivating research in stochastic 
travelling waves are stochastic neural field equations~\cite{Bressloff,Bressloff3,
FaugerasInglis,KuehnRiedler}
\be
\label{eq:nfwaves}
\txtd u = \left[-\alpha u + \int_\cI f(u(t,y))w(\cdot,y)~\txtd y\right]~\txtd t + \sigma g(u)~\txtd W,\quad 
\alpha>0,~u=u(t,x),~\cI\subseteq \R,
\ee
where $w$ is a kernel modelling the connections of the neurons, and other variants of 
neural field equations place the nonlinearity $f$ outside of the integral $f(\int_\cI\ldots)$. 
It has been proven
that~\eqref{eq:nfwaves} has many analogies to classical local (S)PDEs~\cite{KuehnToelle,LaingTroy}.
Furthermore, travelling waves have been studied, particularly in the bistable case, in quite
some detail for stochastic neural fields, see 
e.g.~\cite{BressloffWebber,BressloffWilkerson,InglisMacLaurin,KilpatrickErmentrout,KruegerStannat,Lang}.

An important topic directly related to the bistable setting is the 
case $a=1/2$ for $f=f_3$, so that the PDE has a standing wave. It is well-known
that if we consider a weak diffusion
\be
\label{eq:ACstanding}
\partial_t =\varepsilon^2\partial_x^2 u + u(1-u)(u-1/2)
\ee
then~\eqref{eq:ACstanding} develops quickly, for quite large sets of initial data, several
sharp interfaces of width $\cO(\varepsilon)$ between $0$ and $1$. These interfaces then move 
exponentially slowly on an (approximating) invariant manifold at speed $\cO(\txte^{-K/\varepsilon})$ 
for some constant $K>0$; see e.g.~\cite{CarrPego,KuehnBook1}. Of course, one may then 
ask, how these interfaces form and move in the case, when~\eqref{eq:ACstanding} is perturbed 
by noise. This case has been studied in quite some detail showing that the interfaces 
still form and move~\cite{Funaki,Funaki1,Lee2}. It is anticipated that an invariant manifold 
description still exists~\cite{AntonopoulouBatesBloemkerKarali,AntonopoulouBloemkerKarali},
and that dynamics on this manifold is essentially Brownian motion~\cite{BrassescoDeMasiPresutti}
under suitable conditions.

A very important future direction for research will be to connect waves/patterns for 
SPDEs more closely to applications~\cite{GarciaOjalvoSancho}. Stochastic wave-like
structures in SPDEs have already appeared in an extremely diverse set of modelling contexts
such as neuroscience~\cite{SauerStannat2,Tuckwell1,Tuckwell2}, spin 
glasses~\cite{BrunetDerridaMuellerMunier2}, biological invasions~\cite{Snyder}, predator-prey 
systems~\cite{KhainLinSander,SieberMalchowPetrovskii}, directed polymers~\cite{BrunetDerrida1},
evolutionary biology~\cite{CohenKesslerLevine,Hallatschek}, and 
epidemics~\cite{Linetal,Warrenetal}. Although small fluctuations in the modelling context are 
sometimes just neglected, this is generally a false hypothesis \emph{near instability}, e.g., when 
we are close to propagation failure or when the deterministic PDE part undergoes a bifurcation. We
remark that in this context the precise formulation of the SPDE via modelling will be crucial. The 
bifurcation (or phase/critical transition) aspect has been recognized early on as a key concept in 
SPDEs~\cite{HohenbergHalperin} for steady-state-like patterns. A recently more detailed 
mathematical theory has begun to develop~\cite{GowdaKuehn,KuehnRomano} to use critical slowing down
in combination with stochastic perturbations as early-warning signs for transitions of 
steady-state-like stochastic dynamics. However, a similar idea also seems to have emerged 
early on in the context of numerical simulations of waves in noisy systems~\cite{SchloeglBerry}. 
Therefore, we conjecture that one may very efficiently compare data and SPDE models 
via fluctuation analysis near wave/pattern-forming instabilities. 

\medskip
\textbf{Acknowledgments:} I would like to thank the VolkswagenStiftung for 
support via a Lichtenberg Professorship. I also acknowledge the very helpful comments
of two anonymous referees, of the editor and of Christian Hamster, who thereby helped to 
improve the presentation of the manuscript.

%%%%%%%%%%%%%%%%%%%%%%%%%%%%%%%%%%%%%%%%%%%%%%%%%%%%%%%%%%%%%%%%%%%%%%%%%%%%%%%%%%%%%%%%%%%%%%
{
\scriptsize
\bibliographystyle{plain}
\bibliography{../my_refs}
}
\end{document}